\documentclass[12pt,pdflatex]{sn-jnl}
\usepackage{graphicx}%
\usepackage{subcaption}%
\usepackage{multirow}%
\usepackage{amsmath,amssymb,amsfonts}%
\usepackage{amsthm}%
\usepackage{mathrsfs}%
\usepackage[title]{appendix}%
\usepackage{xcolor}%
\usepackage{textcomp}%
\usepackage{manyfoot}%
\usepackage{booktabs}%
\usepackage{algorithm}%
\usepackage{algorithmicx}%
\usepackage{algpseudocode}%
\usepackage{listings}%
\usepackage{enumitem}
\usepackage{makecell}
\usepackage{caption}

\theoremstyle{thmstyleone}%
\newtheorem*{theorem*}{Theorem}
\newtheorem{theorem}{Theorem}[section]

\newtheorem{lemma}[theorem]{Lemma}% meant for sectionwise numbers
\newtheorem{proposition}[theorem]{Proposition}% 

\theoremstyle{thmstyletwo}%
\newtheorem{example}[theorem]{Example}%
\newtheorem{remark}[theorem]{Remark}%

\theoremstyle{thmstylethree}%

\begin{document}

\title[Damped Traub's method]{Dynamics of damped Traub's method applied to polynomials with two distinct roots}

%%=============================================================%%
%% GivenName	-> \fnm{Joergen W.}
%% Particle	-> \spfx{van der} -> surname prefix
%% FamilyName	-> \sur{Ploeg}
%% Suffix	-> \sfx{IV}
%% \author*[1,2]{\fnm{Joergen W.} \spfx{van der} \sur{Ploeg} 
%%  \sfx{IV}}\email{iauthor@gmail.com}
%%=============================================================%%

\author[1]{\fnm{Subham} \sur{Chatterjee}}\email{subham.acad.res@gmail.com}
\equalcont{These authors contributed equally to this work.}

\author[1]{\fnm{Gorachand} \sur{Chakraborty}}\email{gorachand.chakraborty@skbu.ac.in}
\equalcont{These authors contributed equally to this work.}

\author*[2]{\fnm{Subhasis} \sur{Ghora}}\email{subhasisghora@iiitdm.ac.in}
\equalcont{These authors contributed equally to this work.}

\affil[1]{\orgdiv{Department of Mathematics}, \orgname{Sidho-Kanho-Birsha University}, \orgaddress{\city{Purulia}, \postcode{723104}, \state{West Bengal}, \country{India}}}

\affil[2]{\orgdiv{Department of Sciences and Humanities}, \orgname{Indian Institute of Information
Technology Design and Manufacturing Kancheepuram}, \orgaddress{\city{Chennai}, \postcode{600127}, \state{Tamil Nadu}, \country{India}}}

%\affil[3]{\orgdiv{Department}, \orgname{Organization}, \orgaddress{\street{Street}, \city{City}, \postcode{610101}, \state{State}, \country{Country}}}

%%==================================%%
%% Sample for unstructured abstract %%
%%==================================%%

\abstract{
We study the dynamics of damped Traub's method $T_{p,\delta}(z)=N_p(z)-\delta\frac{p(N_p(z))}{p^\prime(z)}$, one-parameter family that interpolates between Newton's method $N_p$ ($\delta=0$) and Traub's method ($\delta=1$) applied to polynomials $p$ with two distinct roots. Due to the scaling property, it is sufficient to study its dynamics applied to $p_{\vphantom{1}_{n}}(z)=(z+1)^k(z-1)^{n-k}.$ It is proven that for $T_{p_{\vphantom{1}_{2k}},\delta}$, the immediate basins of attraction of $\pm 1$ are completely contained  in the left half-plane and right half-plane, respectively and that $T_{p_{\vphantom{1}_{2k}},\delta}$ is semiconjugate via $z\mapsto z^2$ to a rational map of simpler form. The main focus of the paper is to study the dynamics of $T_{p_{\vphantom{1}_{2}},\delta}$ with real $\delta$. We show that the Julia set of $T_{p_{\vphantom{1}_{2}},\delta}$  is disconnected for $\delta <-4$ and $\delta >12$, while it is connected for some suitable values of $\delta$. The Fatou set is also completely described for $\delta\in(-\infty,-\frac{4}{3}]\cup(0,\frac{3+\sqrt{13}}{2}]\cup[\frac{19+5\sqrt{17}}{8},\infty)$. In particular, it is proved that $T_{p_{\vphantom{1}_{2}},\delta}$ has no Siegel disk or Herman ring  for real parameter values. On the other hand, for certain complex parameters, the existence of invariant Siegel disks for $T_{p_{\vphantom{1}_{2}},\delta}$ has been established. Numerical experiments show a period-doubling in the range $\delta\in(4,\frac{19+5\sqrt{17}}{8})$. It is established that unlike Newton’s method, damped Traub’s method can have extraneous attracting or parabolic cycles, Siegel disks, and disconnected Julia sets even for quadratics.
}


\maketitle

\section{Introduction}\label{sec1}
Nonlinear equations arise in many areas of science such as physics, biology, engineering, and applied mathematics etc. However, finding solutions of nonlinear equations has been a challenge to scientists. In many cases, exact algebraic solutions are nearly impossible to find. To compute approximate solutions, mathematicians have developed iterative root-finding methods. Understanding the dynamics of these methods is essential to identifying stable regions that lead to convergence to the roots.

A root-finding method is a function that associates a polynomial $p$ with a rational function $R_p$ such that any zero $z_0$ of $p$ is an attracting fixed point of $R_p$, meaning $R_p(z_0)=z_0$ and $\lvert R_p^\prime(z_0)\rvert <1$. This implies that if we choose $z$ sufficiently close to $z_0$, then $\{z_n\}$ converges to $z_0$, where $z_n=R_p^n(z)$ for all $n\in\mathbb{N}$. If $z_0$ is not an attracting fixed point then  $\{z_n\}$ may not converge to $z_0$ which is not desired for a root-finding method. 

The most well-known root-finding method is Newton's method, given by $$N_p(z)=z-\frac{p(z)}{p^\prime(z)}$$ for a polynomial $p$. Newton's method has an order of convergence $2$ for simple roots. For higher-order convergence, Traub's method was introduced by J. F. Traub \cite{traub1982iterative}. For a polynomial $p$, Traub's method is given by 
$$
T_p(z)=N_p(z)-\frac{p(N_p(z))}{p^\prime(z)}.
$$
This method has an order of convergence $3$ for simple roots. To bridge these two classical methods into a single continuous family,  damped Traub's method is introduced. For a polynomial $p$, the damped Traub's method is given by 
$$
T_{p,\delta}(z)=N_p(z)-\delta\frac{p(N_p(z))}{p^\prime(z)},
$$ where $\delta$ is a complex number known as the damping factor. For $\delta=0$ and $\delta=1$, the damped Traub's method becomes Newton's method and Traub's method respectively.
In this work, we study the dynamics of damped Traub's method applied to polynomials with two distinct roots. For $\delta\neq 1$, the method has an order of convergence $2$. Even though the order of convergence decreases for $\delta\neq 1$, the regions of convergence can be larger for certain values of $\delta$. Traub’s method offers impressive local third-order convergence, but its performance degrades quickly when initial guesses lie far from a root. The primary advantage of the damped Traub’s method over standard Traub's method is that varying the damping factor provides direct control over fixed points and unstable periodic orbits.

For a rational function $f:\widehat{\mathbb{C}}\to \widehat{\mathbb{C}}$ of degree at least $2$, we consider the sequence of functions $\displaystyle\left\{f^n\right\}_{n\geq 0}$, where $f^n=f^{n-1}\circ f$ and $f^0$ is the identity function. Based on the behaviour of this sequence, we divide the Riemann sphere into two disjoint subsets, namely the Fatou set and the Julia set. The Fatou set of $f$ is defined as $F(f)=\{z\in\widehat{\mathbb{C}}: z\text{ has a neighbourhood } U\text{ such that the family } \displaystyle\left\{f^n\right\} \text{ is normal on } U\}.$
The Julia set is the complement of the Fatou set. Maximal connected subsets of the Fatou set  are known as Fatou components.
For rational functions, if $U$ is a Fatou component, then $U_n=U_m$ for some $n,m\in\mathbb{N}\cup\{0\}\text{ and }n\neq m$, where $U_k$ is the Fatou component containing $f^k(U)$ \cite[Theorem 8.1.2]{beardon2000iteration}. Fatou components satisfying this property are known as pre-periodic Fatou components. 
If for a Fatou component $U$, there exists $p\in \mathbb{N}$ such that $U_p=U$, then we call it a periodic Fatou component. Periodic Fatou components are classified into $4$ types, namely immediate attracting basin, immediate parabolic basin, Siegel disk, and Herman ring.  A Fatou component $U$ of $f$ is said to be completely invariant if $f^{-1}(U)= U=f(U)$. Periodic Fatou components are closely associated with periodic points. A point $z_0\in\mathbb{C}$ is said to be a $p$-periodic point of $f$ if $f^p(z_0)=z_0$ and $p$ is the least natural number satisfying this property. For a $p$-periodic point $z_0$, the value $\lambda=(f^p)^\prime(z_0)$ is known as the multiplier at $z_0$. Periodic points are classified into $3$ categories, namely attracting, indifferent, and repelling, based on the modulus of the multiplier.  
$1$-periodic points are known as fixed points. More details about periodic Fatou components and periodic points are given in Section \ref{S:s2}.
For an attracting $p$-periodic point $z_0$, the attracting basin of $z_0$ is defined as $\mathcal A(z_0)=\{z\in\widehat{\mathbb{C}}:f^{np}(z)\to z_0 \text{ as } n\to\infty\}$. Similarly, for a rationally indifferent $p$-periodic point  $z_0$, the parabolic basin of $z_0$ is defined as $\mathcal P(z_0)=\{z\in\widehat{\mathbb{C}}:f^{np}(z)\to z_0\text{ as } n\to\infty\}\setminus J(f)$. Attracting basins or parabolic basins need not be connected. The component of the attracting basin $\mathcal A(z_0)$ containing the attracting periodic point $z_0$ is $\mathcal A_{z_0}$, the immediate attracting basin of $z_0$. Similarly, the component of the parabolic basin $\mathcal P(z_0)$ whose boundary contains the rationally indifferent periodic point $z_0$ is $\mathcal P_{z_0}$, the immediate parabolic basin of $z_0$. For a root-finding method $R_p$, any fixed point that is not a root of $p$ is known as an extraneous fixed point.  
 If the extraneous fixed point is attracting, then it has a corresponding attracting basin. Thus, choosing the initial value from that basin gives us convergence to a point that is not a root of the parent function. Due to this, root-finding methods with extraneous fixed points have very complicated dynamics. We have shown in Section \ref{S:s3} that damped Traub's method can have extraneous fixed points.
Traub's method and its variations have been studied from a dynamical point of view by many mathematicians. It is shown that for $r_d(z)=z(z^d-1)$, the immediate attracting basins of $T_{r_d}$ corresponding to the roots of $r_d$ are unbounded and simply connected \cite{canela2026boundedness}. In \cite{canela2023basins}, it is proved that for any quadratic polynomial $p$, the immediate basins of attraction of $T_{p,\delta}$ corresponding to the roots of $p$ are simply connected and unbounded for $\delta\in(0,1]$. The order of convergence of damped Traub's method and the stability of the fixed points have been studied in \cite{cordero2016damped}.
 
The Scaling Theorem is a very important aspect of the dynamics of root-finding methods. It allows us to study the dynamics of a class representative to get an idea about the dynamics of the whole class. 
Many well-known root-finding methods, including Halley's method \cite{liu2026dynamics}, Chebyshev's method \cite{ghora2026chebyshev,nayak2025symmetry}, K\"onig's method \cite{buff2003konig}, etc., satisfy the Scaling Theorem \cite{amat2004review}, and the dynamics of these methods have been studied by many researchers. On the other hand, there are some root-finding methods that do not satisfy the Scaling Theorem like Steffensen's method \cite{mandal2026dynamics1,mandal2026dynamics2}. 
Damped Traub's method also satisfies the Scaling Theorem as shown in Lemma \ref{L:scaling}.
In this work, we take $p_{\vphantom{1}_{n}}$ as the polynomial 
$
p_{\vphantom{1}_{n}}(z)=(z+1)^k(z-1)^{n-k}, n,k\in \mathbb{N} \text{ and } n>k.
$ Due to the scaling property of damped Traub's method, for polynomials with exactly two distinct roots, it is sufficient to study the dynamics of $T_{p_{\vphantom{1}_{n}},\delta}$.

The primary focus of this work is to provide a dynamical description of $T_{p_{\vphantom{1}_{2}},\delta}$ for $\delta\in\mathbb{R}$. While most of the other works on damped Traub's method have focused on the connectivity and boundedness of the immediate attracting basins of the roots, our focus has been to describe the Julia set in greater detail. The existence of a disconnected Julia set is quite rare in the dynamics of root-finding methods. Newton's method, when applied to polynomials or transcendental entire functions, always has connected Julia set \cite{baranski2018connectivity}. Halley's method can have disconnected Julia sets for certain transcendental entire functions \cite{cumsille2022disconnected}.
The existing literature on damped Traub's method has primarily focused on the immediate basins of the roots. Canela et al. \cite{canela2023basins} established that these basins are simply connected and unbounded for quadratic polynomials when $\delta\in(0,1]$, while Cordero et al. \cite{cordero2016damped} investigated the order of convergence and the stability of the fixed points. The present work extends these investigations by studying the global dynamics of the associated rational map and providing a description of its Fatou and Julia sets over the real parameter line. The main contributions of this work are summarised as follows.
\begin{enumerate}
    \item[(i)] We establish that the Julia set of $T_{p_{\vphantom{1}_{2}},\delta}$ is disconnected for $\delta<-4$ and $\delta>12$. This provides polynomial examples of a root-finding method of this type whose Julia sets are disconnected, in contrast to Newton's method, for which the Julia set is connected \cite{baranski2018connectivity}.
    \item[(ii)] We prove that $T_{p_{\vphantom{1}_{2}},\delta}$ has no Siegel disks or Herman rings for real values of $\delta$. In addition, for suitable non-real values of the parameter, we construct two invariant Siegel disks.
    \item[(iii)] For $\delta\in(0,\frac{3+\sqrt{13}}{2}]$, we prove that every Fatou component is simply connected, thereby establishing the connectedness of the Julia set. The proof exploits the fact that the imaginary axis is contained in $J(T_{p_2,\delta})$ and separates the two free critical points. This further extends the work presented in \cite{canela2023basins} where this was shown in the range $\delta\in(0,1]$.
    \item[(iv)] For $\delta>\frac{19+5\sqrt{17}}{8}$, we reduce the dynamics of the free critical orbits to a real odd map with a pole. Further, period doubling is investigated numerically in the range $\delta\in(4,\frac{19+5\sqrt{17}}{8})$. A complete analytic classification of the Fatou set for $\delta\in(4,\frac{19+5\sqrt{17}}{8})$ remains open.
\end{enumerate}
In the following theorems, we describe the dynamics of $T_{p_{\vphantom{1}_{2}},\delta}$ in different ranges of $\delta\in\mathbb{R}$.

\begin{theorem}\label{T:(inf,-4)}
Let $\delta\in(-\infty,-4]$.
\begin{enumerate}[label=\alph*)]
    \item If $\delta\in(-\infty,-4)$, then the Fatou set of $T_{p_{\vphantom{1}_{2}},\delta}$ is the union of $\mathcal A(-1)$, $\mathcal A(1)$, and $\mathcal A(\infty)$. Further, the Julia set of $T_{p_{\vphantom{1}_{2}},\delta}$ is disconnected for all $\delta\in(-\infty,-4)$.
    \item  If $\delta=-4$, then the Fatou set is the union of $\mathcal A(-1)$, $\mathcal A(1)$, and $\mathcal P(\infty)$. Also, the Julia set of $T_{p_{\vphantom{1}_{2}},-4}$ is connected.
\end{enumerate}
\end{theorem}
For $\delta \in (-4, -\frac{4}{3}]$, the dynamics is very dependent on the extraneous fixed points. The two extraneous fixed points of $T_{p_{\vphantom{1}_{2}},\delta}$ are denoted by $e_1$ and $e_2$.
\begin{theorem}\label{T:JuliaTp2conn}
    Let $\delta\in (-4,-\frac{4}{3}]$.
    \begin{enumerate}[label=\alph*)]
        \item Then the Julia set of $T_{p_{\vphantom{1}_{2}},\delta}$ is connected.
        \item If $\delta\in (-4,-\frac{4}{3})$. then the Fatou set of $T_{p_{\vphantom{1}_{2}},\delta}$ is the union of $\mathcal A(-1)$, $\mathcal A(1)$, $\mathcal A(e_1)$, and $\mathcal A(e_2)$.
        \item If $\delta=-\frac{4}{3}$, then the Fatou set of $T_{p_{\vphantom{1}_{2}},\delta}$ is the union of $\mathcal A(-1)$, $\mathcal A(1)$, $\mathcal P(e_1)$, and $\mathcal P(e_2)$.
    \end{enumerate}
\end{theorem}
For $\delta\in(-\frac{4}{3},0)$, all fixed points other than $1$ and $-1$ are repelling. Despite this, we show that the Julia set is connected.
\begin{theorem}\label{T:(-4/3,0)}
    Let $\delta\in(-\frac{4}{3},0)$. Then both $\mathcal A_{-1}$ and $\mathcal A_1$ are simply connected and unbounded. Further, the Julia set of $T_{p_{\vphantom{1}_{2}},\delta}$ is connected.
\end{theorem}
A critical point of $R_p$ that is neither a root of $p$ nor a pole of $R_p$, is known as a free critical point. Existence of periodic Fatou component other than the basins of attraction of the roots ensure the existence of free critical points. For suitable values of $\delta$, it is shown that $T_{p_{\vphantom{1}_{2}},\delta}$ has two free critical points and both of them are real and hence, their orbits are also restricted to the real line. The following theorem discusses it.
\begin{theorem}\label{T:0,4}
Let $\delta\in(0,4)$.
\begin{enumerate}[label=\alph*)]
    \item Then the imaginary axis denoted by $\{iy:y\in\mathbb{R}\}$, is contained in the Julia set for $\delta\in (0,4)$.
    \item If $\delta\in(0,1]$, then the Fatou set of $T_{p_{\vphantom{1}_{2}},\delta}$ is the union of $\mathcal A(-1)$ and $\mathcal A(1)$, and the Julia set of $T_{p_{\vphantom{1}_{2}},\delta}$ is connected. Further, both $\mathcal A_{-1}$ and $\mathcal A_1$ are simply connected and unbounded.
    \item If $\delta\in(1,\frac{3+\sqrt{13}}{2}]$, then the Fatou set is the union of $\mathcal A(-1)$ and $\mathcal A(1)$. Also, the Julia set is connected.
\end{enumerate}    
\end{theorem}
In Theorem \ref{T:0,4}, the case for $\delta\in(0,1]$ has previously been discussed in \cite{canela2023basins}, although it used different approach.

\begin{theorem}\label{T:4.8,12}
    Let $\delta\in [\frac{19+5\sqrt{17}}{8},12)$.
    \begin{enumerate}[label=\alph*)]
        \item If $\delta\in(\frac{19+5\sqrt{17}}{8},12)$, then the Fatou set of $T_{p_{\vphantom{1}_{2}},\delta}$ is the union of $\mathcal A(-1),\mathcal A(1)$, and a $2$-cycle of attracting basins.
        \item If $\delta=\frac{19+5\sqrt{17}}{8}$, then the Fatou set is the union of $\mathcal A(-1),\mathcal A(1)$, and a $2$-cycle of parabolic basins.
        \item The Julia set of $T_{p_{\vphantom{1}_{2}},\delta}$ is connected.
    \end{enumerate}
\end{theorem}
For $\delta\in[12,\infty)$, the point at infinity becomes an attracting or rationally indifferent
fixed point.
\begin{theorem}\label{T:12,inf}
Let $\delta\in[12,\infty)$.
\begin{enumerate}[label=\alph*)]
    \item If $\delta\in (12,\infty)$, then the Fatou set of $T_{p_{\vphantom{1}_{2}},\delta}$ is the union of $\mathcal A(-1),\mathcal A(1)$, and $\mathcal A(\infty)$. Further, the Julia set of $T_{p_{\vphantom{1}_{2}},\delta}$ is disconnected.
    \item  If $\delta=12$, then the Fatou set of $T_{p_{\vphantom{1}_{2}},\delta}$ is the union of $\mathcal A(-1),\mathcal A(1)$, and $\mathcal P(\infty)$. Also, the Julia set of $T_{p_{\vphantom{1}_{2}},12}$ is connected.
\end{enumerate} 
\end{theorem}

In Section \ref{S:s2}, we describe some general properties of damped Traub's method applied to polynomials. Section \ref{S:s3} is entirely devoted to studying the dynamics of damped Traub's method applied to $p_{\vphantom{1}_{n}}$, 
primarily focusing on the dynamics of $T_{p_{\vphantom{1}_{2}},\delta}$. 
Throughout this paper, we assume $-D=\{-z:z\in D\}$ for any subset $D$ of $\mathbb{C}$.

\section{Preliminaries}\label{S:s2}
We first describe some basic concepts of complex dynamics. For a periodic point $z_0\in\mathbb{C}$ with multiplier $\lambda$,
if $\lvert \lambda\rvert<1$, then $z_0$ is attracting. If $\lvert \lambda\rvert=1$, then $z_0$ is indifferent. In this case, there are two possibilities. If $\lambda=e^{2\pi i\theta}$ with $\theta\in[0,1)\cap \mathbb{Q}$, then $z_0$ is rationally indifferent. Otherwise, it is irrationally indifferent. If $\lvert \lambda\rvert>1$, then $z_0$ is repelling. The point at infinity is a $p$-periodic point of $f$ if $g^p(0)=0$, where $g(z)=\frac{1}{f(\frac{1}{z})}$, and in that case the multiplier at infinity is given by $(g^p)^\prime(0)$.
For the rational function $f$, the zeros of $f^\prime$ and the multiple poles of $f$ are known as critical points, and the images of critical points are known as critical values.

As mentioned before, periodic Fatou components can be of $4$ types for a rational function. Let $U$ be a $p$-periodic Fatou component of $f$.
\begin{itemize}
    \item If $U$ contains an attracting $p$-periodic point $a$ such that $f^{np}(z)\to a$ as $n\to\infty$ for all $z\in U$, then $U$ is said to be an immediate attracting basin. The immediate attracting basin containing $a$ is denoted by $\mathcal A_a$.
    \item If $\partial U$, the boundary of $U$, contains a rationally indifferent periodic point $r$ such that $f^{np}(z)\to r$ as $n\to\infty$ for all $z\in U$, then $U$ is said to be an immediate parabolic basin. The immediate parabolic basin whose boundary contains $r$ is denoted by $\mathcal P_r$.
    \item If there exists an analytic homeomorphism $\phi:U\to\mathbb{D}$ such that $\phi\circ f^p\circ \phi^{-1}(z)=e^{2\pi i\theta}z$ for all $z\in \mathbb{D}$ and some $\theta\in[0,1)\cap \mathbb{Q}^c$, then $U$ is known as a Siegel disk.
    \item If $\mathbb{A}=\{z:1<\lvert z\rvert<r\}, r>1$, and there exists an analytic homeomorphism $\phi:U\to\mathbb{A}$ such that $\phi\circ f^p\circ \phi^{-1}(z)=e^{2\pi i\theta}z$ for all $z\in \mathbb{A}$ and some $\theta\in[0,1)\cap \mathbb{Q}^c$, then $U$ is known as a Herman ring.
\end{itemize}

 For a connected domain $U$, the connectivity of $U$, denoted by $c(U)$, is the number of components of $\widehat{\mathbb{C}}\setminus U$. To study the connectivity of Fatou components, the Riemann-Hurwitz formula is particularly useful.
\begin{theorem*}[Riemann-Hurwitz formula \cite{canela2026boundedness}]
    Let $U,V$ be two connected domains, and let $f:U\to V$ be a degree-$d$ proper map branched over $k$ critical points, counting multiplicities. Then
    $$c(U)-2=d(c(V)-2)+k.$$
\end{theorem*}

We show that damped Traub's method applied to polynomials satisfies the scaling property.
\begin{lemma}\label{L:scaling}
    Let $p$ be a polynomial of degree at least $2$ and $h(z)=az+b$, with $a,b\in\mathbb{C}$, and $a\neq 0$. Then for any $\lambda\neq 0$, $(h^{-1}\circ T_{p,\delta}\circ h) \equiv T_{q,\delta}$ where $q(z)=\lambda(p\circ h)(z)$ for all $z\in\mathbb{C}$.
\end{lemma}
\begin{proof}
    Note that $T_{q,\delta}(z)=N_q(z)-\delta\frac{q(N_q(z))}{q^\prime(z)}=z-\frac{q(z)}{q^\prime (z)}-\delta\frac{q(N_q(z))}{q^\prime(z)}$. Since $q(z)=\lambda(p\circ h)(z)$, we have $$N_q(z)=z-\frac{q(z)}{q^\prime(z)}=z-\frac{p(h(z))}{p^\prime(h(z))h^\prime(z)}=z-\frac{p(h(z))}{ap^\prime(h(z))}.$$ Also, 
    $$\begin{aligned}
        q(N_q(z))& =\lambda p(h(N_q(z)))=\lambda p\left(az+b-\frac{p(h(z))}{p^\prime(h(z))}\right)\\
        &=\lambda p\left(h(z)-\frac{p(h(z))}{p^\prime(h(z))}\right)=\lambda p(N_p(h(z)).
    \end{aligned}    
    $$ This implies that $T_{q,\delta}(z)=z-\frac{p(h(z))}{ap^\prime(h(z))}-\delta \frac{p(N_p(h(z))}{ap^\prime(h(z))}$. It follows that $(h\circ T_{q,\delta})(z)=N_p(h(z))-\delta\frac{p(N_p(h(z))}{p^\prime(h(z))}=(T_{p,\delta}\circ h)(z)$. Therefore, $$(h^{-1}\circ T_{p,\delta}\circ h) (z)= T_{q,\delta}(z).$$
\end{proof}
\begin{remark}
    To simplify the dynamics of a function, the scaling property is very useful. For example, for any polynomial $p$ of degree $n$, the damped Traub's method $T_{p,\delta}$ is conjugate to   $T_{q,\delta}$ for some monic polynomial $q$ of degree $n$.
\end{remark}
    In this work, we study the dynamics of $T_{p,\delta}$ for a monic polynomial $p$ with exactly two distinct roots.  We assume $p$ is a monic polynomial of degree $n$ with two distinct roots $\alpha$ and $\beta$ such that $\alpha$ is a root of multiplicity $k$ with $ 0 < k < n$. The polynomial $p$ is of the form $p(z)=(z-\alpha)^k(z-\beta)^{n-k}$. Considering $h(z)=\frac{\beta-\alpha}{2}z+\frac{\beta+\alpha}{2}$ and $\lambda=\frac{2^n}{(\beta-\alpha)^n}$, we have $\lambda(p\circ h)(z)=(z+1)^k(z-1)^{n-k}$. Without loss of generality, we can assume $p(z)=(z+1)^k(z-1)^{n-k}$.
    If a polynomial $q$ is of the form $q(z)=(z+1)^{n-k}(z-1)^k$, then using scaling, it can be shown that $q$ and $p$ are affinely conjugate. Throughout the article, we consider such polynomials $p_{\vphantom{1}_{n}}$ with the following form  $$p_{\vphantom{1}_{n}}(z)=(z+1)^k(z-1)^{n-k}.$$

Since the damped Traub's method reduces to Newton's method for $\delta=0$, we assume $\delta\neq 0$ throughout the paper.

The damped Traub's method applied to ${p_{\vphantom{1}_{n}}}$ is given by
\begin{equation}\label{E:TP}
    T_{p_{\vphantom{1}_{n}},\delta}(z)=N_{p_{\vphantom{1}_{n}}}(z)-\delta\frac{p(N_{p_{\vphantom{1}_{n}}}(z))}{p^\prime(z)}=z-B(z)C(z)
\end{equation} where 
\begin{equation*}
    B(z)=\frac{z^2-1}{L(z)^{n+1}}
\end{equation*}
  and 
  \begin{equation*}
      C(z)=L(z)^n+\delta(L(z)-z+1)^k(L(z)-z-1)^{n-k}
  \end{equation*}  
  for 
  \begin{equation*}
    L(z)=nz+(n-2k).  
  \end{equation*}
Clearly, the finite fixed points of $T_{p_{\vphantom{1}_{n}},\delta}$ are $\pm 1$ and the solutions of $C(z)=0$. The nature of the fixed points $-1$ and $1$ is described in the following lemma.
\begin{lemma}\label{L:multiplier-11}
    The multipliers at the fixed points $-1$ and $1$ of  $T_{p_{\vphantom{1}_{n}},\delta}$ are given by $\frac{k-1}{k}-\delta\frac{(k-1)^k}{k^{k+1}}$ and $\frac{n-k-1}{n-k}-\delta\frac{(n-k-1)^{n-k}}{(n-k)^{n-k+1}}$, respectively.
\end{lemma}
\begin{proof}
    From Eqn. \eqref{E:TP}, we get
\begin{equation}\label{E:derTP}
    T_{p_{\vphantom{1}_{n}},\delta}^\prime(z)=1-B(z)C^\prime(z)-B^\prime(z)C(z).
\end{equation}
Note that 
\begin{equation*}
   B^\prime(z)=\frac{2zL(z)-n(n+1)(z^2-1)}{L(z)^{n+2}} 
\end{equation*}
and 
\begin{equation*}
   C^\prime(z)=n^2L(z)^{n-1}+\delta(n-1)(L(z)-z+1)^{k-1}(L(z)-z-1)^{n-k-1}\big(n(L(z)-z)+n-2k\big). 
\end{equation*}
From Eqn. \eqref{E:derTP}, we have $T_{p_{\vphantom{1}_{n}},\delta}^\prime(-1)=1-B^\prime(-1)C(-1)$ and $T_{p_{\vphantom{1}_{n}},\delta}^\prime(1)=1-B^\prime(1)C(1)$. The rest of the proof follows from simple calculations.
\end{proof}
\begin{remark}\label{R:k1super}
    It follows from Lemma \ref{L:multiplier-11} that if $k=1$ then $-1$ is always superattracting fixed point of $T_{p_{\vphantom{1}_{n}},\delta}$. Similarly, if $n-k=1$, then $1$ is a superattracting fixed point of $T_{p_{\vphantom{1}_{n}},\delta}$. 
\end{remark}
 
We have the following results from \cite[Lemma 2.1]{canela2023basins}.
 
\begin{lemma}\label{L:range-1}
\begin{enumerate}[label=\alph*)]
    \item For $k\neq 1$, the fixed point $-1$ of $T_{p_{\vphantom{1}_{n}},\delta}$ is attracting, rationally indifferent, and repelling if $-\frac{k^k}{(k-1)^k}<\delta<(2k-1)\frac{k^k}{(k-1)^k}$, $\delta=-\frac{k^k}{(k-1)^k}$ or $\delta=(2k-1)\frac{k^k}{(k-1)^k}$, and $\delta<-\frac{k^k}{(k-1)^k}$ or $\delta>(2k-1)\frac{k^k}{(k-1)^k}$ respectively.
    \item For $n-k\neq 1$, the fixed point $1$ of $T_{p_{\vphantom{1}_{n}},\delta}$ is attracting, rationally indifferent, and repelling if $-\frac{(n-k)^{n-k}}{(n-k-1)^{n-k}}<\delta<(2(n-k)-1)\frac{(n-k)^{n-k}}{(n-k-1)^{n-k}}$, $\delta=-\frac{(n-k)^{n-k}}{(n-k-1)^{n-k}}$ or $\delta=(2(n-k)-1)\frac{(n-k)^{n-k}}{(n-k-1)^{n-k}}$ and $\delta<-\frac{(n-k)^{n-k}}{(n-k-1)^{n-k}}$ or $\delta>(2(n-k)-1)\frac{(n-k)^{n-k}}{(n-k-1)^{n-k}}$ respectively.
\end{enumerate}
\end{lemma}
Since $T_{p_{\vphantom{1}_{n}},\delta}$ is a rational function, it is possible for the point at infinity to be fixed point of $T_{p_{\vphantom{1}_{n}},\delta}$. In fact, the point at infinity is always a fixed point of $T_{p_{\vphantom{1}_{n}},\delta}$ for $\delta\neq \frac{n^n}{(n-1)^{n-1}}$. For $\delta= \frac{n^n}{(n-1)^{n-1}}$, it no longer remains a fixed point. We have the following result from \cite[Lemma 2.2]{canela2023basins}.
 
\begin{lemma}\label{L:rangeinfty}
    Let $\delta\in\mathbb{R}\setminus\{\frac{n^n}{(n-1)^{n-1}},0\}$. 
    \begin{enumerate}[label=\alph*)]
        \item The point at infinity is attracting, rationally indifferent, and repelling fixed point of $T_{p_{\vphantom{1}_{n}},\delta}$, if $\delta\in(-\infty,-\frac{n^n}{(n-1)^{n}})\cup((2n-1)\frac{n^n}{(n-1)^{n}},\infty)$, $\delta=-\frac{n^n}{(n-1)^{n}}$ or $\delta=(2n-1)\frac{n^n}{(n-1)^{n}}$, and $\delta\in(-\frac{n^n}{(n-1)^{n}},(2n-1)\frac{n^n}{(n-1)^{n}})$ respectively.
    \end{enumerate}
\end{lemma}
\begin{remark}\label{R:rangeinfty}
    It follows from Lemma \ref{L:rangeinfty}, that for $n=2$, the point at infinity is attracting if $\delta\in (-\infty,-4)\cup (12,\infty)$, rationally indifferent if $\delta=-4$ or $\delta=12$ and repelling if $\delta\in(-4,12)\setminus\{0,4\}$. Recall that for $\delta=0$, the method reduces to the Newton's method and for $\delta=4$, the point at infinity is not a fixed point of $T_{p_{\vphantom{1}_{n}},\delta}$.
\end{remark}
The following result can be proved using \cite[Lemma 3.4]{nayak2025symmetry} and  \cite[Lemma 4.4]{ghora2026chebyshev}.
 
\begin{lemma}\label{L:general1}
    Let $R$ be a rational map of degree at least $2$ such that the point at infinity is either a repelling or a rationally indifferent fixed point of $R$. Let $U$ be an invariant immediate attracting basin or immediate parabolic basin of $R$ such that $\infty\in \partial U$. Then $\partial U$ contains at least one pole of $R$.
\end{lemma}
The following result, proved in \cite[Lemma 3.5]{nayak2025symmetry} states about the connectedness of the Julia set of a rational map under certain condition.
 
\begin{lemma}\label{L:general2}
    Let $R$ be a rational map of degree at least $2$ such that the point at infinity is either a repelling or a rationally indifferent fixed point of $R$. If an unbounded component of the Julia set contains all poles of $R$, then the Julia set is connected.
\end{lemma}

\section{Dynamics of $T_{p_{\vphantom{1}_{n}},\delta}$}\label{S:s3}
%We start this section by proving that the Fatou set of $T_{p,\delta}$ is entire complex plane then $p$ has some particular form. 
%\begin{proposition}
%    Let $p$ be a non-constant polynomial. If $F(T_{p,\delta})=\mathbb{C}$, then $T_{p,\delta}$ is conjugate to $T_{q,\delta}$ for some $q(z)=z^n$ for $n\geq 2$.
%\end{proposition}
%\begin{proof}
%    Since $\mathbb{C}$ is
%the only Fatou component, it contains at most one periodic cycle concluding that $T_{p,\delta}$ has
%at most one finite periodic point, and it is fixed. Any finite point is not a pole of $T_{p,\delta}$ as  $J(T_{p,\delta})=\{\infty\}$. Hence $T_{p,\delta}$ is a polynomial. If the degree of $T_{p,\delta}$ is at
%least $2$, then $\infty$ is a superattracting fixed point contained in the Fatou set,
%contradicting $J(T_{p,\delta})=\{\infty\}$. Thus $T_{p,\delta}$ is linear. If $p$ has $k(\geq 2)$ many distinct roots, then $T_{p,\delta}$ has at least two finite fixed points. This is not possible. This implies that $p$ has exactly one zero, say $z_0$. So it must be of the form $p(z)=c(z-z_0)^n$ for some $c\in\mathbb{C}\setminus\{0\}$ and $n\in\mathbb{N}$. 
%Note that $n\geq 2$. Otherwise, the map $T_{p,\delta}(z)=z_0$ for all $z\in\mathbb{C}$. This implies that $F(T_{p,\delta})=\widehat{\mathbb{C}}$ which is a contradiction. Considering $h(z)=z+z_0$, we get $q(z)=\frac{1}{c}(p\circ h)(z)=z^n$ for some $n\geq 2$. By scaling property $T_{p,\delta}$ is conjugate to $T_{q,\delta}$ under the map $h$.
%\end{proof}

To describe the behaviour of the fixed points of $T_{p_{\vphantom{1}_{n}},\delta}$ where $p_{\vphantom{1}_{n}}(z)=(z+1)^k(z-1)^{n-k}$, the following lemma is useful. 
 
\begin{lemma}\label{L:increasing}
    Let $g_1(x)=-\frac{x^x}{(x-1)^x}$ and $g_2(x)=(2x-1)\frac{x^x}{(x-1)^x}$ for $x\geq 2$. Then both $g_1$ and $g_2$ are strictly increasing functions.
\end{lemma}
\begin{proof}
    Considering $r(x)=\frac{x^x}{(x-1)^x}$, we get $\frac{r^\prime(x)}{r(x)}=\log{(1+\frac{1}{x-1})}-\frac{1}{x-1}$. Hence, $$g_1^\prime(x)=r(x)\left(\frac{1}{x-1}-\log{\left(1+\frac{1}{x-1}\right)}\right)>0$$ for all $x\geq 2$. Now, for $g_2$ taking logarithm and differentiating, we get 
    $$
    \begin{aligned}
        g_2^\prime(x)&=g_2(x)\left(\frac{2}{2x-1}+\frac{r^\prime(x)}{r(x)}\right)\\
        & =g_2(x)\left(\frac{2}{2x-1}+\log{(1+\frac{1}{x-1})}-\frac{1}{x-1}\right)\\
        &=g_2(x)\left(\log{(1+\frac{1}{x-1})}-\frac{\frac{1}{(x-1)^2}}{2+\frac{1}{x-1}}\right).
    \end{aligned}
    $$ 
    Let $t=\frac{1}{x-1}$ which implies that $0<t\leq 1$. Then $\log(1+t)-\frac{t^2}{2+t}>\frac{2t}{2+t}-\frac{t^2}{2+t}=\frac{t(2-t)}{2+t}>0$. Consequently, we have $g_2^\prime(x)>0$ for all $x\geq 2$. Hence, both $g_1$ and $g_2$ are strictly increasing function for all $x\geq 2$.
\end{proof}
 Since the zeros of $p_{\vphantom{1}_{n}}$, $-1$ and $1$, are attracting only in particular ranges of $\delta$, it is necessary to obtain the common range for which both zeros are attracting. For this the next lemma is important. 
\begin{lemma}\label{L:attracatingrange1}
    Let $2k<n$ and $k,n-k\neq 1$. Then the following are true for $T_{p_{n},\delta}$.
    \begin{enumerate}[label=\alph*)]
        \item If $\delta<-\frac{k^k}{(k-1)^k}$ or $\delta>(2(n-k)-1)\frac{(n-k)^{n-k}}{(n-k-1)^{n-k}}$, then both $-1$ and $1$ are repelling.
        \item If $\delta=-\frac{k^k}{(k-1)^k}$, then $-1$ is rationally indifferent and $1$ is repelling.
        \item If $-\frac{k^k}{(k-1)^k}<\delta<-\frac{(n-k)^{n-k}}{(n-k-1)^{n-k}}$, then $-1$ is attracting and $1$ is repelling.
        \item If $\delta=-\frac{(n-k)^{n-k}}{(n-k-1)^{n-k}}$, then $-1$ is attracting and $1$ is rationally indifferent.
        \item If $-\frac{(n-k)^{n-k}}{(n-k-1)^{n-k}}<\delta<(2k-1)\frac{k^k}{(k-1)^k}$, then both $-1$ and $1$ are attracting.
        \item If $\delta=(2k-1)\frac{k^k}{(k-1)^k}$, then $-1$ is rationally indifferent and $1$ is attracting.
        \item If $(2k-1)\frac{k^k}{(k-1)^k}<\delta<(2(n-k)-1)\frac{(n-k)^{n-k}}{(n-k-1)^{n-k}}$, then $-1$ is repelling and $1$ is attracting.
        \item If $\delta=(2(n-k)-1)\frac{(n-k)^{n-k}}{(n-k-1)^{n-k}}$, then $-1$ is repelling and $1$ is rationally indifferent.
    \end{enumerate}
\end{lemma}
\begin{proof}
    Since $2k<n$, it follows that $k<n-k$. Following Lemma \ref{L:increasing} we have $g_1(k)<g_1(n-k)<0<g_2(k)<g_2(n-k)$.
    \begin{enumerate}[label=\alph*)]
        \item If $\delta<-\frac{k^k}{(k-1)^k}$, then $\delta<g_1(k)<g_1(n-k)$. If $\delta>(2(n-k)-1)\frac{(n-k)^{n-k}}{(n-k-1)^{n-k}}$, then $\delta>g_2(n-k)>g_2(k)$. So, by Lemma \ref{L:range-1}, both $-1$ and $1$ are repelling.
        \item Considering $\delta=-\frac{k^k}{(k-1)^k}=g_1(k)<g_1(n-k)$, the fixed points $-1$ and $1$ are rationally indifferent and repelling respectively.
        \item When $-\frac{k^k}{(k-1)^k}<\delta<-\frac{(n-k)^{n-k}}{(n-k-1)^{n-k}}$, then $g_1(k)<\delta<g_1(n-k)<0$. So, the fixed points $-1$ and $1$ are attracting and repelling respectively.
        \item If $\delta=-\frac{(n-k)^{n-k}}{(n-k-1)^{n-k}}$, then $g_1(k)<g_1(n-k)=\delta<0$. Hence, the fixed points $-1$ and $1$ are attracting and rationally indifferent respectively. 
        \item Since $-\frac{(n-k)^{n-k}}{(n-k-1)^{n-k}}<\delta<(2k-1)\frac{k^k}{(k-1)^k}$, $g_1(k)<g_1(n-k)<\delta<g_2(k)<g_2(n-k)$. So, both $-1$ and $1$ are attracting.
        \item For $\delta=(2k-1)\frac{k^k}{(k-1)^k}$, $0<\delta=g_2(k)<g_2(n-k)$. Therefore, the fixed points $-1$ and $1$ are rationally indifferent and attracting respectively.
        \item As $(2k-1)\frac{k^k}{(k-1)^k}<\delta<(2(n-k)-1)\frac{(n-k)^{n-k}}{(n-k-1)^{n-k}}$, $g_2(k)<\delta<g_2(n-k)$. Thus, the fixed points $-1$ and $1$ are repelling  and attracting respectively.
        \item For $\delta=(2(n-k)-1)\frac{(n-k)^{n-k}}{(n-k-1)^{n-k}}$, $\delta=g_2(n-k)>g_2(k)$. Hence, the fixed points $-1$ and $1$ are repelling and rationally indifferent respectively.
    \end{enumerate}
\end{proof}
\begin{remark}
    If $2k>n$ and $k,(n-k)\neq 1$, then $q_n(z)=\lambda(p_n\circ h)(z)=(z+1)^{n-k}(z-1)^k$ where $\lambda=\frac{1}{(-1)^n}$ and $h(z)=-z$. Since $n-k<k$, using Lemma \ref{L:attracatingrange1} we can show that both $-1$ and $1$ are attracting fixed points of $q_n$ if $-\frac{k^k}{(k-1)^k}<\delta<(2(n-k)-1)\frac{(n-k)^{n-k}}{(n-k-1)^{n-k}}$. By scaling property, it follows that $T_{p_{\vphantom{1}_{n}},\delta}(z)=h\circ T_{q_{\vphantom{1}_{n}},\delta}\circ h^{-1}(z)=-T_{q_{\vphantom{1}_{n}},\delta}(-z)$, and $T_{p_{\vphantom{1}_{n}},\delta}^\prime(z)=T_{q_{\vphantom{1}_{n}},\delta}^\prime(-z)$. So, it follows that $\lvert T_{p_{\vphantom{1}_{n}},\delta}^\prime(1)\rvert=\lvert T_{q_{\vphantom{1}_{n}},\delta}^\prime(-1)\rvert$ and $\lvert T_{p_{\vphantom{1}_{n}},\delta}^\prime(-1)\rvert=\lvert T_{q_{\vphantom{1}_{n}},\delta}^\prime(1)\rvert$. Therefore, the fixed points $-1$ and $1$ are attracting for $T_{p_{\vphantom{1}_{n}},\delta}$ if $-\frac{k^k}{(k-1)^k}<\delta<(2(n-k)-1)\frac{(n-k)^{n-k}}{(n-k-1)^{n-k}}$. Hence, it is sufficient to consider only the case $n>2k$ described in Lemma \ref{L:attracatingrange1}. 
\end{remark}
%Similarly we have the following result.
 
%\begin{lemma}\label{L:attracatingrange2}
%    Let $2k>n$ and $k,n-k\neq 1$. Then the following are true for $T_{p_{n},\delta}$.
%    \begin{enumerate}[label=\alph*)]
%        \item If $\delta<-\frac{(n-k)^{n-k}}{(n-k-1)^{n-k}}$ or $\delta>(2k-1)\frac{k^k}{(k-1)^k}$, then both $-1$ and $1$ are repelling.
%        \item If $\delta=-\frac{(n-k)^{n-k}}{(n-k-1)^{n-k}}$, then $-1$ is repelling and $1$ is rationally indifferent.
%        \item If $-\frac{(n-k)^{n-k}}{(n-k-1)^{n-k}}<\delta<-\frac{k^k}{(k-1)^k}$, then $-1$ is repelling and $1$ is attracting.
%        \item If $\delta=-\frac{k^k}{(k-1)^k}$, then $-1$ is rationally indifferent and $1$ is attracting.
%        \item If $-\frac{k^k}{(k-1)^k}<\delta<(2(n-k)-1)\frac{(n-k)^{n-k}}{(n-k-1)^{n-k}}$, then both $-1$ and $1$ are attracting.
%        \item If $\delta=(2(n-k)-1)\frac{(n-k)^{n-k}}{(n-k-1)^{n-k}}$, then $-1$ is attracting and $1$ is rationally indifferent.
%        \item If $(2(n-k)-1)\frac{(n-k)^{n-k}}{(n-k-1)^{n-k}}<\delta<(2k-1)\frac{k^k}{(k-1)^k}$, then $-1$ is attracting and $1$ is repelling.
%        \item If $\delta=(2k-1)\frac{k^k}{(k-1)^k}$, then $-1$ is rationally indifferent and $1$ is repelling.
%    \end{enumerate}
%\end{lemma}
From Lemma \ref{L:attracatingrange1}, we can see that $T_{p_{\vphantom{1}_{n}},\delta}$ is suitable as a root-finding method only if $$-\frac{(n-k)^{n-k}}{(n-k-1)^{n-k}}<\delta<(2k-1)\frac{k^k}{(k-1)^k}.$$

If $n=2k$, then from Lemma \ref{L:range-1}, both $-1$ and $1$ are attracting fixed point if $$-\frac{k^k}{(k-1)^k}<\delta<(2k-1)\frac{k^k}{(k-1)^k}.$$

The fixed points of $T_{p_{\vphantom{1}_{n}},\delta}$ which are not the roots of $p_{\vphantom{1}_{n}}$ are the solutions of $C(z)=0$. But due to the extremely broad nature of $C$, it is very difficult to find the exact solutions of $C(z)=0$ for arbitrary choices of $n,k$. However, if we choose $n=2k$, then we can find these finite extraneous fixed points.
 
\begin{lemma}
    Let $n=2k$. Then there exists $2k$ finite extraneous fixed points of $T_{p_{\vphantom{1}_{n}},\delta}$ and they are 
    given by $$z=\pm\displaystyle\sqrt{\frac{\lvert\delta\rvert^{\frac{1}{k}}e^{n_r\pi i}}{\lvert\delta\rvert^{\frac{1}{k}}(2k-1)^2e^{n_r\pi i}-4k^2}},$$ 
    and the multiplier at $z$ is given by 
    $$\displaystyle\frac{\lvert\delta\rvert^{\frac{2}{k}}(2k-1)^2(1-k)e^{2n_r\pi i}+2k(2k-1)^2\lvert\delta\rvert^{\frac{1}{k}}e^{n_r\pi i}-4k^3}{k\lvert\delta\rvert^{\frac{1}{k}}e^{n_r\pi i}},$$
    where 
$n_r = \begin{cases} 
\frac{2r+1}{k} & \text{if } \delta > 0 \\
\frac{2r}{k} & \text{if } \delta<0 
\end{cases}$ for $r= 0,1,2,\dots,k-1$.
\end{lemma}
\begin{proof}
    Considering $n=2k$, we have $L(z)=nz$ and $C(z)=n^nz^n+\delta\big((n-1)^2z^2-1\big)^k$. The zeros of $C$ are the solution of 
    \begin{equation}\label{E:2k}
        \left(\frac{z^2}{(2k-1)^2z^2-1}\right)^{k}=\frac{-\delta}{(2k)^{2k}}.
    \end{equation}
    At a zero of $C$, the derivative is  $C^\prime(z)=(2k)^{2k+1}z^{2k-1}\left(1-\frac{(2k-1)^2z^2}{(2k-1)^2z^2-1}\right)$, and hence, we have from Eqn. \eqref{E:derTP}, 
    \begin{equation}\label{E:multiextra}        T_{p_{\vphantom{1}_{n}},\delta}^\prime(z)=1+\frac{z^2-1}{z^2((2k-1)^2z^2-1)}.
    \end{equation} 
    
\textbf{Case I:} Assuming $\delta>0$, it follows from Eqn. \eqref{E:2k} that $\frac{z^2}{(2k-1)^2z^2-1}=\frac{\lvert\delta\rvert^{\frac{1}{k}}}{4k^2}e^{\frac{2r+1}{k} \pi i}$ for $r=0,1,,\dots,k-1$ implying  that $z^2=\frac{\lvert\delta\rvert^{\frac{1}{k}}e^{\frac{2r+1}{k}\pi i}}{\lvert\delta\rvert^{\frac{1}{k}}(2k-1)^2e^{\frac{2r+1}{k}\pi i}-4k^2}$. Thus, the finite extraneous fixed points are given by $$z=\pm\displaystyle\sqrt{\frac{\lvert\delta\rvert^{\frac{1}{k}}e^{\frac{2r+1}{k}\pi i}}{\lvert\delta\rvert^{\frac{1}{k}}(2k-1)^2e^{\frac{2r+1}{k}\pi i}-4k^2}}  \text{ for } r=0,1,2,\dots,k-1.$$ Putting $z^2=\frac{\lvert\delta\rvert^{\frac{1}{k}}e^{\frac{2r+1}{k}\pi i}}{\lvert\delta\rvert^{\frac{1}{k}}(2k-1)^2e^{\frac{2r+1}{k}\pi i}-4k^2}$ in Eqn. \eqref{E:multiextra} and simplifying, we get $$T_{p_{\vphantom{1}_{n}},\delta}^\prime(z)=\displaystyle\frac{\lvert\delta\rvert^{\frac{2}{k}}(2k-1)^2(1-k)e^{\frac{4r+2}{k}\pi i}+2k(2k-1)^2\lvert\delta\rvert^{\frac{1}{k}}e^{\frac{2r+1}{k}\pi i}-4k^3}{k\lvert\delta\rvert^{\frac{1}{k}}e^{\frac{2r+1}{k}\pi i}}.$$

\textbf{Case II:} Considering $\delta<0$, it follows from Eqn. \eqref{E:2k} that $\frac{z^2}{(2k-1)^2z^2-1}=\frac{\lvert\delta\rvert^{\frac{1}{k}}e^{\frac{2r}{k}\pi i}}{4k^2}$ for $r=0,1,2,\dots,k-1.$ This implies that $z^2=\frac{\lvert\delta\rvert^{\frac{1}{k}}e^{\frac{2r}{k}\pi i}}{\lvert\delta\rvert^{\frac{1}{k}}(2k-1)^2e^{\frac{2r}{k}\pi i}-4k^2}.$  The finite extraneous fixed points are given by  $$z=\pm\displaystyle\sqrt{\frac{\lvert\delta\rvert^{\frac{1}{k}}e^{\frac{2r}{k}\pi i}}{\lvert\delta\rvert^{\frac{1}{k}}(2k-1)^2e^{\frac{2r}{k}\pi i}-4k^2}} \text{ for } r=0,1,2,\dots,k-1.$$ Putting $z^2=\frac{\lvert\delta\rvert^{\frac{1}{k}}e^{\frac{2r}{k}\pi i}}{\lvert\delta\rvert^{\frac{1}{k}}(2k-1)^2e^{\frac{2r}{k}\pi i}-4k^2}$ in Eqn. \eqref{E:multiextra} we get $$T_{p_{\vphantom{1}_{n}},\delta}^\prime(z)=\displaystyle\frac{\lvert\delta\rvert^{\frac{2}{k}}(2k-1)^2(1-k)e^{\frac{4r}{k}\pi i}+2k(2k-1)^2\lvert\delta\rvert^{\frac{1}{k}}e^{\frac{2r}{k}\pi i}-4k^3}{k\lvert\delta\rvert^{\frac{1}{k}}e^{\frac{2r}{k}\pi i}}.$$   
\end{proof}
\begin{example}
    If we take $k=2$ and $\delta=-3$, then $n_r=0,1$ and the extraneous fixed points are $z=\pm\sqrt{\frac{\sqrt{3}}{16-9\sqrt{3}}}i\approx\pm 2.0515i$ and $z=\pm\sqrt{\frac{\sqrt{3}}{16+9\sqrt{3}}}\approx\pm 0.2342$. At the extraneous fixed points $z=\pm 0.2342$, the multiplier is given by $M\approx35.03$. So, the extraneous fixed points $z=\pm 0.2342$ are repelling. However, at the extraneous fixed points $z=\pm 2.0515i$, the multiplier is given by $M\approx0.968$. So, the extraneous fixed points $z=\pm 2.0515i$ are attracting.
\end{example}
 Let $\mathbb{H}^-=\{x+iy:x<0, y\in\mathbb{R}\}$ and $\mathbb{H}^+=\{x+iy:x>0, y\in\mathbb{R}\}$. 
\begin{lemma}\label{L:Tpninvline}
    If $n=2k$, $k>1$, and $-\frac{k^k}{(k-1)^k}<\delta<(2k-1)\frac{k^k}{(k-1)^k}$, then $\mathcal A_{-1}\subset \mathbb{H}^-$ and $\mathcal A_1\subset \mathbb{H}^+$.
\end{lemma}
\begin{proof}
    As $n=2k$, we have $$T_{p_{\vphantom{1}_{2k}},\delta}(z)=z-\frac{z^2-1}{(2k)^{2k+1}z^{2k+1}}\Big((2k)^{2k}z^{2k}+\delta\big((2k-1)^2z^2-1\big)^k\Big).$$
 %   Thus, $T_{p_{\vphantom{1}_{2k}},\delta}(-z)=-T_{p_{\vphantom{1}_{2k}},\delta}(z)$.
    Since $n=2k$ and $-\frac{k^k}{(k-1)^k}<\delta<(2k-1)\frac{k^k}{(k-1)^k}$, both $-1$ and $1$ are attracting fixed points.
    Also, for all $iy$ such that $y\neq 0$,
    $T_{p_{\vphantom{1}_{2k}},\delta}(iy)=ig_k(y)$ where $$g_k(y)=y-\frac{y^2+1}{(2k)^{2k+1}y^{2k+1}}\Big((2k)^{2k}y^{2k}+\delta\big((2k-1)^2y^2+1\big)^k\Big).$$
    This implies that $T_{p_{\vphantom{1}_{2k}},\delta}\Big(\big\{iy:y\in\mathbb{R}\setminus\{0\}\big\}\Big)\subset \big\{iy:y\in\mathbb{R}\big\}$. If $\mathcal A_{-1}$ contains a point in $\mathbb{H}^+$, then there exists some point $iy_0\in \mathcal A_{-1}$. But the forward orbit of $iy_0$ is contained in the imaginary axis. This is a contradiction. So, $\mathcal A_{-1}\subset \mathbb{H}^-$. Similarly, if $\mathcal A_{1}$ contains a point in $\mathbb{H}^-$, then there must be some point $iy_1\in\mathcal A_1$ which is not possible implying $\mathcal A_1\subset \mathbb{H}^+$.
\end{proof}
The conclusion of Lemma \ref{L:Tpninvline} holds even if $k=1$ and $n=2$, as shown in Lemma \ref{T:Tp2civ}.
 
\begin{proposition}\label{P:semiconj}
    If $n=2k$, then $T_{p_{\vphantom{1}_{2k}},\delta}$ is semiconjugate to $$R(z)=z\left(1-\frac{z-1}{(2k)^{2k+1}z^{k+1}}\Big((2k)^{2k}z^{k}+\delta\big((2k-1)^2z-1\big)^k\Big)\right)^2$$ under the map $\phi(z)=z^2$.
\end{proposition}
\begin{proof}
    We have $$
    \begin{aligned}
        \phi\circ T_{p_{\vphantom{1}_{2k}},\delta}(z)& =(T_{p_{\vphantom{1}_{2k}},\delta})^2\\
     &= z^2+\frac{(z^2-1)^2}{(2k)^{4k+2}z^{4k+2}}\left((2k)^{2k}z^{2k}+\delta\left((2k-1)^2z^2-1\right)^k\right)^2 \\
     & - \frac{2(z^2-1)}{(2k)^{2k+1}z^{2k}}\left((2k)^{2k}z^{2k}+\delta\left((2k-1)^2z^2-1\right)^k\right)\\
 & = z^2\Big(1+\frac{(z^2-1)^2}{(2k)^{4k+2}z^{4k+4}}\big((2k)^{2k}z^{2k}+\delta((2k-1)^2z^2-1)^k\big)^2 \\
 & -\frac{2(z^2-1)}{(2k)^{2k+1}z^{2k+2}}\big((2k)^{2k}z^{2k}+\delta((2k-1)^2z^2-1)^k\big)\Big)\\
 &= R\circ \phi(z).
    \end{aligned}   
    $$
\end{proof}
For large values of $k$, it might be helpful to study the dynamics of $R$ as described in Proposition \ref{P:semiconj}, rather than of $T_{p_{\vphantom{1}_{2k}},\delta}$ directly. After studying the dynamics for $R$, one can use semiconjugacy relations to transfer the dynamical properties back to the dynamics of $T_{p_{\vphantom{1}_{2k}},\delta}$.
\subsection{Dynamics of $T_{p_{\vphantom{1}_{2}},\delta}$}
For the next part, we assume $n=2$. Since $n=2$, the value of $k$ must be $1$, and $p_{\vphantom{1}_{2}}(z)=z^2-1$.
Therefore,
\begin{equation}\label{E:Tp2}
    T_{p_{\vphantom{1}_{2}},\delta}(z)=z-\frac{z^2-1}{8z^3}\big((4+\delta)z^2-\delta\big).
\end{equation}
$T_{p_{\vphantom{1}_{2}},\delta}$ is a degree $4$ rational function for all $\delta\neq 4$ and has $6$ critical points counting multiplicities for all $\delta\neq 4$. $T_{p_{\vphantom{1}_{2}},4}$ is a degree $3$ rational function and has $4$ critical points counting multiplicities. From Eqn. \eqref{E:Tp2}, the finite fixed points of $T_{p_{\vphantom{1}_{2}},\delta}$ are $\pm 1$ and $\pm\sqrt{\frac{\delta}{4+\delta}}$ for $\delta\neq -4$. For $\delta=-4$, the only finite fixed points of $T_{p_{\vphantom{1}_{2}},\delta}$ are $\pm 1$. The roots of $p_{\vphantom{1}_{2}}$, the points $\pm 1$ are always superattracting by Remark \ref{R:k1super}. Also, for $\delta\neq 4$, the point at infinity is also a fixed point of $T_{p_{\vphantom{1}_{2}},\delta}$ by Remark \ref{R:rangeinfty}. The point at infinity is an attracting fixed point for $\delta\in(-\infty,-4)\cup(12,\infty)$ by Remark \ref{R:rangeinfty}. 

We also have,
\begin{equation}\label{E:derTp2}
    T_{p_{\vphantom{1}_{2}},\delta}^\prime(z)=\frac{z^2-1}{8z^4}\big((4-\delta)z^2-3\delta\big).
\end{equation}
From Eqn. \eqref{E:derTp2}, it follows that $T_{p_{\vphantom{1}_{2}},\delta}$ has two free critical points $$c_1=\sqrt{\frac{3\delta}{4-\delta}} \text{ and } c_2=-\sqrt{\frac{3\delta}{4-\delta}}$$ for $\delta\neq 1,4$. If $\delta=4$, then $T_{p_{\vphantom{1}_{2}},\delta}$ has no free critical point. If $\delta=1$, then $c_1$ and $c_2$ coincide with $1$ and $-1$ respectively. Since $0$ is a multiple pole of $T_{p_{\vphantom{1}_{2}},\delta}$, it is also a critical point.
\begin{remark}
    As discussed earlier $T_{p_{\vphantom{1}_{2}},4}$ has no free critical points. However, in this case $\{0,\infty\}$ forms a superattracting $2$-cycle. Since $T_{p_{\vphantom{1}_{2}},4}$ has three distinct critical points ($1,-1$ with multiplicity one and $0$ with multiplicity two), in this case the only periodic Fatou components are $\mathcal A_{-1}$ and $\mathcal A_1$ and the $2$-periodic superattracting cycle of Fatou components $\{\mathcal A_0,\mathcal A_\infty\}$ where $\mathcal A_0$ contains $0$ and $\mathcal A_\infty$ contains $\infty$.  Therefore, $$F(T_{p_{\vphantom{1}_{2}},4})=\mathcal A(-1)\cup\mathcal A(1)\cup \mathcal{A}(0)\cup \mathcal{A}(\infty).$$ Also, the immediate basins $\mathcal A_{-1}$ and $\mathcal A_1$ are simply connected since they contain exactly one critical point each \cite{campos2020connectivity}. For the map $T_{p_{\vphantom{1}_{2}},4}^2:\mathcal A_0\to \mathcal A_0$, the immediate basin $\mathcal A_0$ contains only one critical point $0$ of $T_{p_{\vphantom{1}_{2}},4}^2$. Therefore by B\"{o}ttcher's theorem $\mathcal A_0$ is simply connected \cite[Theorem 9.3]{milnor2006dynamics}. Now, let $c(\mathcal A_\infty)=s$ and the degree of the proper map $T_{p_{\vphantom{1}_{2}},4}:\mathcal A_\infty\to \mathcal A_0$ be $d^\prime$. By Riemann-Hurwitz formula we have $s-2=-d^\prime$. Simplifying this we get $s+d^\prime=2$. Since, $s$ and $d^\prime$ are positive integers, we have $s=1$ and $d^\prime=1$. Hence, the basins $\mathcal A_0$ and $ \mathcal A_\infty$ are also simply connected. As rest of the Fatou components are preperiodic with no critical points, every Fatou component is simply connected using Riemann-Hurwitz formula. This implies that $J(T_{p_{\vphantom{1}_{2}},4})$ is connected.
\end{remark}
We now describe the nature of the two extraneous fixed points. 
 
\begin{lemma}\label{L:Tp2extra}
The finite extraneous fixed points of $T_{p_{\vphantom{1}_{2}},\delta}$ are attracting, rationally indifferent or repelling if $\delta\in (-4,-\frac{4}{3})$, $\delta=-\frac{4}{3}$ or $\delta\in(-\infty,-4)\cup(-\frac{4}{3},\infty)$ respectively.
%\begin{enumerate}[label=\alph*)]
%    \item attracting if $\delta\in (-4,-\frac{4}{3})$,
%    \item rationally indifferent if $\delta=-\frac{4}{3}$, and
%    \item repelling if $\delta\in(-\infty,-4)\cup(-\frac{4}{3},\infty)$.
%\end{enumerate}
\end{lemma}
\begin{proof}
    The finite extraneous fixed points of $T_{p_{\vphantom{1}_{2}},\delta}$ are given by $\pm\sqrt{\frac{\delta}{4+\delta}}$. From Eqn. \eqref{E:derTp2}, the multiplier at both extraneous fixed points is $2+\frac{4}{\delta}$. So if $\delta\in (-4,-\frac{4}{3})$, then $1>2+\frac{4}{\delta}>-1$. Hence, the extraneous fixed points are attracting.

    If $\delta=-\frac{4}{3}$, then $2+\frac{4}{\delta}=-1$. Thus, both extraneous fixed points are rationally indifferent.

    For $\delta\in(-\infty,-4)\cup(0,\infty)$ or $\delta\in (-\frac{4}{3},0)$, $2+\frac{4}{\delta}>1$ and $2+\frac{4}{\delta}<-1$, respectively. Therefore, both extraneous fixed points are repelling.
\end{proof}
We denote the two finite extraneous fixed points as $e_1$ and $e_2$ such that $$e_1=\sqrt{\frac{\delta}{4+\delta}} \text{ and }e_2=-\sqrt{\frac{\delta}{4+\delta}}.$$
 
\begin{lemma}\label{L:Tp2symmetry}
   $T_{p_{\vphantom{1}_{2}},\delta}^n(-z)=-T_{p_{\vphantom{1}_{2}},\delta}^n(z)$ and $T_{p_{\vphantom{1}_{2}},\delta}^n(\overline{z})=\overline{T_{p_{\vphantom{1}_{2}},\delta}^n(z)}$ for all $n\in\mathbb{N}$ and $\delta\in\mathbb{R}\setminus\{0\}$.
\end{lemma}
\begin{proof}
    As $\delta\in\mathbb{R}\setminus\{0\}$, it follows that $T_{p_{\vphantom{1}_{2}},\delta}(-z)=-T_{p_{\vphantom{1}_{2}},\delta}(z)$ and $T_{p_{\vphantom{1}_{2}},\delta}(\overline{z})=\overline{T_{p_{\vphantom{1}_{2}},\delta}(z)}$ for all $z\in\mathbb{C}$. Now, by the Principle of Mathematical Induction, we have $T_{p_{\vphantom{1}_{2}},\delta}^n(-z)=-T_{p_{\vphantom{1}_{2}},\delta}^n(z)$ and $T_{p_{\vphantom{1}_{2}},\delta}^n(\overline{z})=\overline{T_{p_{\vphantom{1}_{2}},\delta}^n(z)}$ for all $n\in\mathbb{N}$.
\end{proof}
 
\begin{lemma}\label{L:Tp2invline}
$T_{p_{\vphantom{1}_{2}},\delta}^n(iy)=ig^n(y)$ for all $n\in\mathbb{N}$, where
\begin{equation}\label{E:g}
   g(y)=y-\frac{(y^2+1)\big((4+\delta)y^2+\delta\big)}{8y^3} 
\end{equation}
for all $y\in\mathbb{R}\setminus\{0\}$. 
\end{lemma}
\begin{proof}
    We have $T_{p_{\vphantom{1}_{2}},\delta}(iy)=i\left(y-\frac{(y^2+1)\big((4+\delta)y^2+\delta\big)}{8y^3}\right)$. Thus, $T_{p_{\vphantom{1}_{2}},\delta}(iy)=ig(y)$. Inductively, it follows that $T_{p_{\vphantom{1}_{2}},\delta}^n(iy)=ig^n(y)$ for all $y\in\mathbb{R}\setminus\{0\}$. 
\end{proof}
%\begin{lemma}
%    $T_{p_{\vphantom{1}_{2}},\delta}$ is semiconjugate to the map $R(z)=z\left( 1- \frac{(z-1)((4+\delta)z-\delta)}{8z^2}\right)^2$ under $\phi(z)=z^2$.
%\end{lemma}
%\begin{proof}
%    We have $R\circ \phi(z)=R(z^2)=z^2-\frac{(z^2-1)((4+\delta)z^2-\delta)}{4z^2}+\frac{(z^2-1)^2((4+\delta)z^2-\delta)^2}{64z^6}=(T_{p_{\vphantom{1}_{2}},\delta}(z))^2=\phi\circ T_{p_{\vphantom{1}_{2}},\delta}(z)$.
%\end{proof}
 
\begin{proposition}\label{T:noSDHR}
The Fatou set of $T_{p_{\vphantom{1}_{2}},\delta}$ contains no Siegel disk or Herman ring for any $\delta\in\mathbb{R}\setminus\{0\}$.
\end{proposition}
\begin{proof}
    Let us assume that $T_{p_{\vphantom{1}_{2}},\delta}$ has a Siegel disk $S$. Then $\partial S$ is contained in the closure of the forward orbit of some critical point. However, $-1\in \mathcal A_{-1}$, $1\in \mathcal A_1$ and $T_{p_{\vphantom{1}_{2}},\delta}(0)=\infty$. So, it is sufficient to check the orbits of the free critical points.

    If $\delta\in(-\infty,0)\cup(4,\infty)$, then both free critical points lie in $\{iy:y\neq 0\}$. Therefore, by Lemma \ref{L:Tp2invline}, the forward orbits of the free critical points also lie in $\{iy:y\in\mathbb{R}\}$. Thus, $\partial S\subset i\mathbb{R}$. This is only possible if either $S=\mathbb{H}^-$ or $S=\mathbb{H}^+$. But $-1\in \mathbb{H}^-$ and $1\in \mathbb{H}^+$. This is a contradiction. So, for $\delta\in(-\infty,0)\cup(4,\infty)$, $T_{p_{\vphantom{1}_{2}},\delta}$ has no Siegel disk.

    If $\delta\in(0,4)\setminus\{1\}$, then both free critical points are real. As $T_{p_{\vphantom{1}_{2}},\delta}(\mathbb{R}\setminus\{0\})\subset \mathbb{R}$, the forward orbits of both free critical points lie in $\mathbb{R}$. Thus $\partial S=\mathbb{R}$. So, either $S=\{x+iy:y<0\}$ or $S=\{x+iy:y>0\}$. However, $\mathcal A_1\cap \{x+iy:y<0\}\neq \emptyset$ and $\mathcal A_1\cap \{x+iy:y>0\}\neq \emptyset$. This is a contradiction. So, for $\delta\in(0,4)$, $T_{p_{\vphantom{1}_{2}},\delta}$ has no Siegel disk. Since  $T_{p_{\vphantom{1}_{2}},\delta}$ has no free critical points for $\delta=1$ or $4$, $T_{p_{\vphantom{1}_{2}},\delta}$ has no Siegel disk. 
    
    Similarly, $T_{p_{\vphantom{1}_{2}},\delta}$ has no Herman ring.
\end{proof}
We have already proved that Siegel disks do not exist for $T_{p_{\vphantom{1}_{2}},\delta}$ for real parameter values. However, if we choose the parameter $\delta$ as non-real complex number, then it is possible for invariant Siegel disks to exist in $F(T_{p_{\vphantom{1}_{2}},\delta})$. Let $\theta=\frac{1}{a_1+\frac{1}{a_2+\frac{1}{a_3+\dots}}}$ be a continued fraction. Then by rational convergents of $\theta$, we mean $\frac{p_n}{q_n}=\frac{1}{a_1+\frac{1}{a_2+\frac{1}{a_3+\frac{1}{\ddots+\frac{1}{a_n}}}}}$. An irrational number $\theta$ is called a Brjuno number when $\Sigma_{n=1}^{\infty} \frac{\log{q_{n+1}}}{q_n}<\infty$. We denote the set of Brjuno numbers by $\mathcal B$.
\begin{proposition}
    Let $\theta\in\mathcal B$ and $\delta=\frac{4}{e^{2\pi i\theta}-2}$. Then $F(T_{p_{\vphantom{1}_{2}},\delta})$ contains two invariant Siegel disks surrounding $e_1$ and $e_2$.
\end{proposition}
\begin{proof}
    We have shown that the multiplier at $e_1$ and $e_2$ is $2+\frac{4}{\delta}$. As $\delta=\frac{4}{e^{2\pi i\theta}-2}$, the multipliers at both extraneous fixed points are given by $e^{2\pi i\theta}$. So, both extraneous fixed points are irrationally indifferent with multiplier $e^{2\pi i\theta}$ where $\theta\in\mathcal B$. Hence, by Brjuno-R\"ussman theorem \cite[Theorem 3.36]{branner2014quasiconformal}, $T_{p_{\vphantom{1}_{2}},\delta}$ is locally linearizable near both $e_1$ and $e_2$. Hence, $F(T_{p_{\vphantom{1}_{2}},\delta})$ contains two invariant Siegel disks surrounding $e_1$ and $e_2$.
\end{proof}
It should be noted that the symmetry of the dynamics over the real line observed for real parameters is not valid for complex parameters. This can be seen from Figure \ref{fig:T_p_2,s}.
\begin{figure}[H]
    \centering
    \includegraphics[trim = 5cm 0cm 5cm 0cm, clip, scale=1]{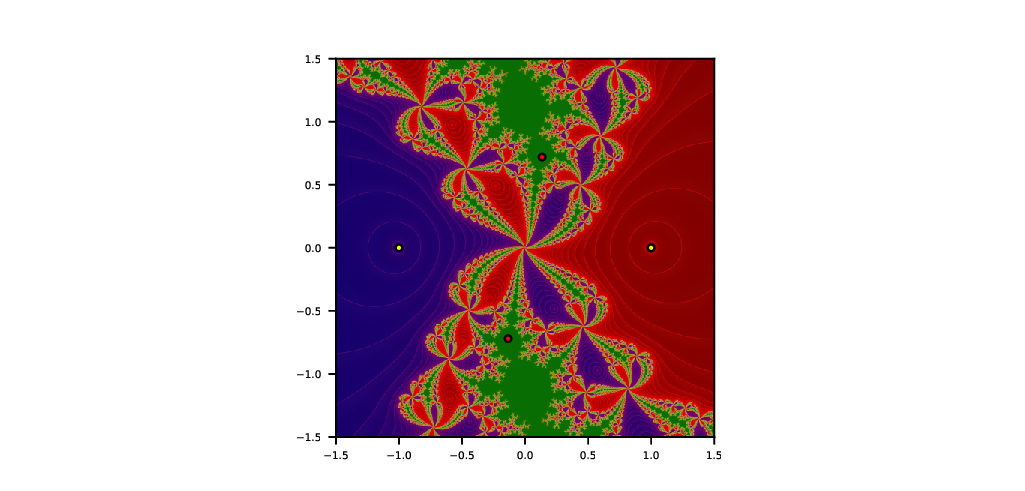}
    \caption{Dynamical plane of $T_{p_{\vphantom{1}_{2}},\delta}$ where $\delta=\frac{4}{e^{(\sqrt{5}-1)\pi i}-2}$. The blue region is $\mathcal A(-1)$ and the red region is $\mathcal A(1)$. The red dots are the two extraneous fixed points $e_1$ and $e_2$, and the green regions surrounding them are the two invariant Siegel disks.}
    \label{fig:T_p_2,s}
\end{figure}

We show that the imaginary axis separates the immediate attracting basins $\mathcal A_{-1}$ and $ \mathcal A_{1}$.
 \begin{lemma}\label{T:Tp2civ}
For $T_{p_{\vphantom{1}_{2}},\delta}$ with $\delta\in\mathbb{R}\setminus\{0\}$, the following are true. 
    \begin{enumerate}[label=\alph*)]
        \item $\mathcal A_{-1}\subset \mathbb{H}^-$ and $\mathcal A_1\subset \mathbb{H}^+$.
        \item $\mathcal A_{-1}$ and $\mathcal A_1$ are not completely invariant.
    \end{enumerate}
\end{lemma}
\begin{proof}
For $p_2$, $T_{p_{\vphantom{1}_{2}},\delta}$ has superattracting fixed points at $z=-1$ and $z=1$ by Remark \ref{R:k1super}.
\begin{enumerate}[label=\alph*)]
    \item  Let us assume that $\mathcal A_{-1}\cap \mathbb{H}^+\neq \emptyset$. Since $\mathcal A_{-1}$ is connected, there must exist some $y_0\in\mathbb{R}\setminus\{0\}$ such that $iy_0\in \mathcal A_{-1}$. But $T_{p_{\vphantom{1}_{2}},\delta}^n(iy_0)\subset \{iy:y\in\mathbb{R}\}$ for all $n\in\mathbb{N}$  by Lemma \ref{L:Tp2invline}. This contradicts $iy_0\in \mathcal A_{-1}$. Therefore $\mathcal A_{-1}\subset \mathbb{H}^-$. Similarly, we can show that $\mathcal A_1\subset \mathbb{H}^+$.
    \item  From Eqn. \eqref{E:Tp2}, we can express $T_{p_{\vphantom{1}_{2}},\delta}$ as
    \begin{equation}\label{E:Tp2inverse-1}
        T_{p_{\vphantom{1}_{2}},\delta}(z)=\frac{(z+1)^2}{8z^3}\big((4-\delta)z^2+2\delta z-\delta\big)-1
    \end{equation}
    and 
    \begin{equation}\label{E:Tp2inverse1}
        T_{p_{\vphantom{1}_{2}},\delta}(z)=\frac{(z-1)^2}{8z^3}\big((4-\delta)z^2-2\delta z-\delta\big)+1.
    \end{equation}
So, for $\delta\neq 4$, from Eqn. \eqref{E:Tp2inverse-1}, $T_{p_{\vphantom{1}_{2}},\delta}(z)=-1$ if either $z=-1$, $z=\frac{\delta}{\delta-4}+2\frac{\sqrt{\delta}}{\delta-4}$ or $z=\frac{\delta}{\delta-4}-2\frac{\sqrt{\delta}}{\delta-4}$. Let $$r_{-1}=\frac{\delta}{\delta-4}+2\frac{\sqrt{\delta}}{\delta-4} \text{ and } r_{-1}^\prime=\frac{\delta}{\delta-4}-2\frac{\sqrt{\delta}}{\delta-4}.$$ 
\textbf{Case 1:} For $\delta>0$, all the preimages of $-1$ are real. If $\delta>4$, then both $r_{-1}$ and $r_{-1}^\prime$ are positive real numbers. So, from a) neither $r_{-1}$ nor $r_{-1}^\prime$ lie in $\mathcal A_{-1}$. If $0<\delta<4$, then $r_{-1}^\prime>0$. So $r_{-1}^\prime\notin \mathcal A_{-1}$. For $\delta=4$, $T_{p_{\vphantom{1}_{2}},4}=-1$ if either $z=-1$ or $z=\frac{1}{2}$. As $\frac{1}{2}>0$, $\mathcal A_{-1}$ does not contain $\frac{1}{2}$.

\textbf{Case 2:} For $\delta<0$, $2\frac{\sqrt{\delta}}{\delta-4}\in \{iy:y\in\mathbb{R}\}$. So, for both $r_{-1}$ and $r_{-1}^\prime$, the real part $\frac{\delta}{\delta-4}>0$ as $\delta<0$. Therefore, all preimages of $-1$ except $-1$ itself lies in $\mathbb{H}^+$ and hence does not lie in $\mathcal A_{-1}$. 

Combining all the cases, $\mathcal A_{-1}$ is not completely invariant for any non-zero real $\delta$. Similarly, using Eqn. \eqref{E:Tp2inverse1}, $\mathcal A_1$ is not completely invariant for all non-zero real $\delta$.

\end{enumerate}
\end{proof}
\begin{remark}\label{R:1pre}
    It is evident from Eqn. \eqref{E:Tp2inverse1}, that $T_{p_{\vphantom{1}_{2}},\delta}^{-1}(1)=\{1, \frac{\delta}{4-\delta}+2\frac{\sqrt{\delta}}{4-\delta},\frac{\delta}{4-\delta}-2\frac{\sqrt{\delta}}{4-\delta}\}.$ We denote $$r_1=\frac{\delta}{4-\delta}+2\frac{\sqrt{\delta}}{4-\delta} \text{ and }r_1^\prime=\frac{\delta}{4-\delta}-2\frac{\sqrt{\delta}}{4-\delta}.$$
\end{remark}
\begin{remark}
For all $\delta\in\mathbb{R}\setminus\{0\}$, there is a component of $T_{p_{\vphantom{1}_{2}},\delta}^{-1}(\mathcal A_{-1})$ that is disjoint from $\mathcal A_{-1}$ due to the previous part. Similarly, there is a component of $T_{p_{\vphantom{1}_{2}},\delta}^{-1}(\mathcal A_{1})$ that is disjoint from $\mathcal A_{1}$. Therefore, $F(T_{p_{\vphantom{1}_{2}},\delta})$ has at least four Fatou components. However, if $J(T_{p_{\vphantom{1}_{2}},\delta})$ is a line, then $F(T_{p_{\vphantom{1}_{2}},\delta})$ have exactly two Fatou components. Hence, $J(T_{p_{\vphantom{1}_{2}},\delta})$ is not a line for any $\delta\in\mathbb{R}\setminus\{0\}$.
\end{remark}
Recall that a repelling fixed point or a rationally indifferent fixed point with multiplier $1$ is known as a weakly repelling fixed point \cite{shishikura2009connectivity}.
\begin{proof}[Proof of Theorem \ref{T:(inf,-4)}]
From Remark \ref{R:rangeinfty}, the point at infinity is attracting if $\delta\in(-\infty,-4)$ and is rationally indifferent with multiplier $1$ if $\delta=-4$.
\begin{enumerate}[label=\alph*)]
    \item From Lemma \ref{L:Tp2invline} we know that $$g(y)=y-\frac{(y^2+1)\big((4+\delta)y^2+\delta\big)}{8y^3}$$ for all $y\in\mathbb{R}\setminus\{0\}$. 
    For $\delta<-4$, $(4+\delta)y^2+\delta<0$ for all $y\in\mathbb{R}\setminus\{0\}$. Since $g(y)-y=-\frac{(y^2+1)\big((4+\delta)y^2+\delta\big)}{8y^3}$, we have $g(y)>y$ for all $y>0$ and $g(y)<y$ for all $y<0$ for $\delta<-4$. Hence, $g$  does not have any real fixed point for $\delta<-4$. Now $$g^\prime(y)=\frac{(4-\delta)y^4+2(2+\delta)y^2+3\delta}{8y^4}.$$ So, $g^\prime(y)=0$ if and only if $y^2=-1$ or $y^2=\frac{3\delta}{\delta-4}$. Since $y\in\mathbb{R}$, $y^2\neq -1$. Therefore, $g^\prime(y)=0$ if and only if either $y=-\sqrt{\frac{3\delta}{\delta-4}}$ or $y=\sqrt{\frac{3\delta}{\delta-4}}$. Let $c=\sqrt{\frac{3\delta}{\delta-4}}$. If $y>c$ or $y<-c$, then $g^\prime(y)>0$. So, for all $y>c$, $g(y)>y$ and $g^\prime(y)>0$. Therefore, for all $y>c$, $g^n(y)\to\infty$ as $n\to\infty$. Similarly, for all $y<-c$, $g^n(y)\to-\infty$ as $n\to\infty$. When $\delta<-4$, then the free critical points of $T_{p_{\vphantom{1}_{2}},\delta}$ are $c_1=ic$ and $c_2=-ic$. As $c>0$, $g(c)>c$. Thus $g^n(c)\to\infty$ as $n\to\infty$. Also, $g^n(-c)\to-\infty$ as $n\to\infty$. Therefore, by Lemma \ref{L:Tp2invline}, $T_{p_{\vphantom{1}_{2}},\delta}^n( c_1)\to \infty$ and $T_{p_{\vphantom{1}_{2}},\delta}^n(c_2)\to \infty$ as $n\to\infty$. 
    %Also, for all $y^2>\frac{3\delta}{\delta-4}$, $T^n_{p_{\vphantom{1}_{2}},\delta}(iy)\to\infty$ as $n\to\infty$.
    So, $\mathcal A(\infty)$ contains both free critical points. Also, $0\in \mathcal A(\infty)$. So, every critical point of $T_{p_{\vphantom{1}_{2}},\delta}$ lies in one of $\mathcal A(-1)$, $\mathcal A(1)$ or $\mathcal A(\infty)$. So, $$F(T_{p_{\vphantom{1}_{2}},\delta})=\mathcal A(-1)\cup \mathcal A(1)\cup \mathcal A(\infty).$$

    To show that the Julia set of $T_{p_{\vphantom{1}_{2}},\delta}$ is disconnected we only need to show that the imaginary axis lies in $\mathcal A_\infty$. We have already proved that for $y>c$ or $y<-c$, $iy\in \mathcal A_\infty$. As $g^\prime(c)=0$ and $g^{\prime\prime}(c)>0$, $g$ has a local minima at $c$ and $g(y)>g(c)$ for all $0<y<c$. As $\delta\in(-\infty,-4)$, $g(c)> c$. So, for all $0<y\leq c$, $g(y)>c$. So, by the previous parts $g^n(y)\to\infty$ as $n\to\infty$ for all $0<y\leq c$. Similarly for $-c\leq y<0$, $g^n(y)\to-\infty$ as $n\to\infty$. So, for all $y\neq 0$, $iy\in \mathcal A_\infty$. Since $0\in \mathcal A(\infty)$ and $iy\in \mathcal A_\infty$ for all $y\neq 0$, $0\in \mathcal A_\infty$. So, the imaginary axis is contained in $\mathcal A_\infty$. However, due to Lemma \ref{T:Tp2civ}, $\mathcal{A}_{-1}\subset \mathbb{H}^-$ and $\mathcal{A}_{1}\subset \mathbb{H}^+$. Consequently, the Julia set of $T_{p_{\vphantom{1}_{2}},\delta}$ is disconnected.    
    \item  For $\delta=-4$, the point at infinity is a rationally indifferent fixed point of $T_{p_{\vphantom{1}_{2}},-4}$ and hence $0\in J(T_{p_{\vphantom{1}_{2}},-4})$. As in the first case, we have $$F(T_{p_{\vphantom{1}_{2}},-4})=\mathcal A(-1)\cup \mathcal A(1)\cup \mathcal P(\infty).$$ $T_{p_{\vphantom{1}_{2}},-4}$ does not have any extraneous fixed point. So, $T_{p_{\vphantom{1}_{2}},-4}$ has only one weakly repelling fixed point at infinity. However, we know from \cite{shishikura2009connectivity} that if a rational function has disconnected Julia set, then it must have at least $2$ distinct weakly repelling fixed points. So, $J(T_{p_{\vphantom{1}_{2}},-4})$ is connected.
\end{enumerate}
\end{proof}
 
\begin{figure}[H]
    \centering
    \includegraphics[trim = 5cm 0cm 5cm 0cm, clip, scale=.8]{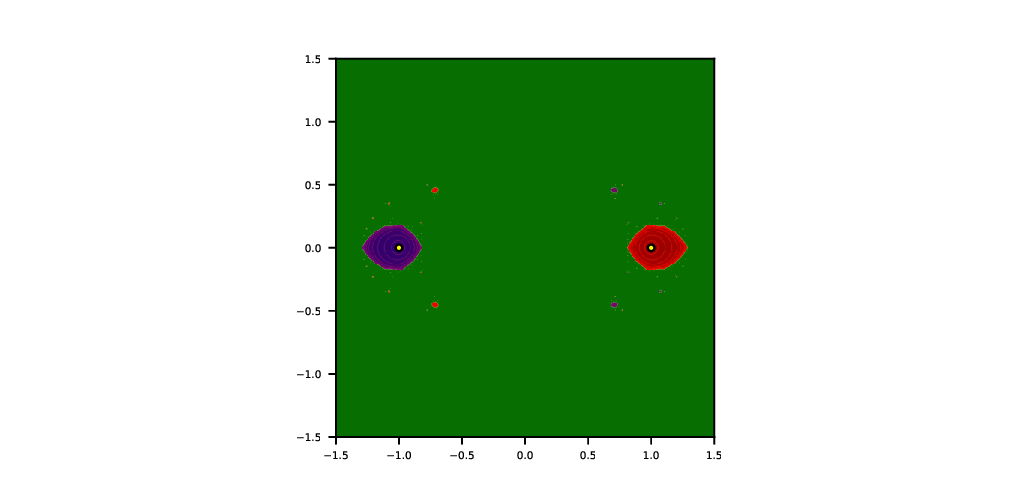}
    \caption{Dynamical plane of $T_{p_{\vphantom{1}_{2}},-10}$. The blue region is $\mathcal A(-1)$ and the red region is $\mathcal A(1)$. The green region is $\mathcal A(\infty)$. The points $-1$ and $1$ are represented by the two yellow dots.}
    \label{fig:T_p_2,-10}
\end{figure}
\vspace{-30pt}
\begin{proof}[Proof of Theorem \ref{T:JuliaTp2conn}]
We have already shown that $-1,1$ are superattracting fixed points of $T_{p_{\vphantom{1}_{2}},\delta}$.
\begin{enumerate}[label=\alph*)]
    \item If $\delta\in (-4,-\frac{4}{3})$, then the two finite extraneous fixed points are attracting by Lemma \ref{L:Tp2extra}. If $\delta=-\frac{4}{3}$, then both extraneous fixed points are rationally indifferent with multiplier $-1$. So, the point at infinity is the only weakly repelling fixed point of $T_{p_{\vphantom{1}_{2}},\delta}$ if $\delta\in (-4,-\frac{4}{3}]$. Therefore, $J(T_{p_{\vphantom{1}_{2}},\delta})$ is connected for all $\delta\in (-4,-\frac{4}{3}]$.
    \item If $\delta\in (-4,-\frac{4}{3})$, then each of $\mathcal A_{-1}$, $\mathcal A_1$, $\mathcal A_{e_1}$ and $\mathcal A_{e_2}$ must contain at least one critical point. Also, the critical point $0$ lies in the Julia set since in this case the point at infinity is a repelling fixed point by Remark \ref{R:rangeinfty}. So, all remaining critical points lie in the immediate attracting basins of the finite fixed points. Hence, only periodic components of $F(T_{p_{\vphantom{1}_{2}},\delta})$ are $\mathcal A_{-1}$, $\mathcal A_1$, $\mathcal A_{e_1}$ and $\mathcal A_{e_2}$. Thus, $$F(T_{p_{\vphantom{1}_{2}},\delta})=\mathcal A(-1)\cup \mathcal A(1)\cup \mathcal A(e_1)\cup \mathcal A(e_2).$$
    \item Similarly, for $\delta=-\frac{4}{3}$, we can show that $$F(T_{p_{\vphantom{1}_{2}},-\frac{4}{3}})=\mathcal A(-1)\cup \mathcal A(1)\cup \mathcal P(e_1)\cup \mathcal P(e_2).$$
\end{enumerate} 
\end{proof}
\begin{figure}[H]
  \centering  
  \begin{subfigure}[b]{0.48\textwidth}
    \centering
    \includegraphics[trim = 5cm 0cm 5cm 0cm, clip, scale=1]{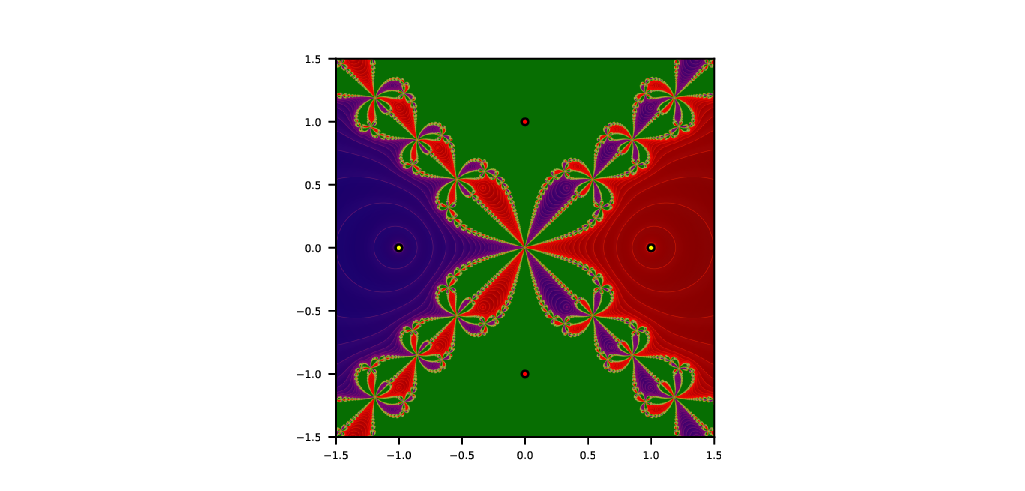}
    \caption{$T_{p_{\vphantom{1}_{2}},-2}$}
    \label{fig:T_p_2,-2}
  \end{subfigure}
  \hfill
  \begin{subfigure}[b]{0.48\textwidth}
    \centering
    \includegraphics[trim = 5cm 0cm 5cm 0cm, clip, scale=1]{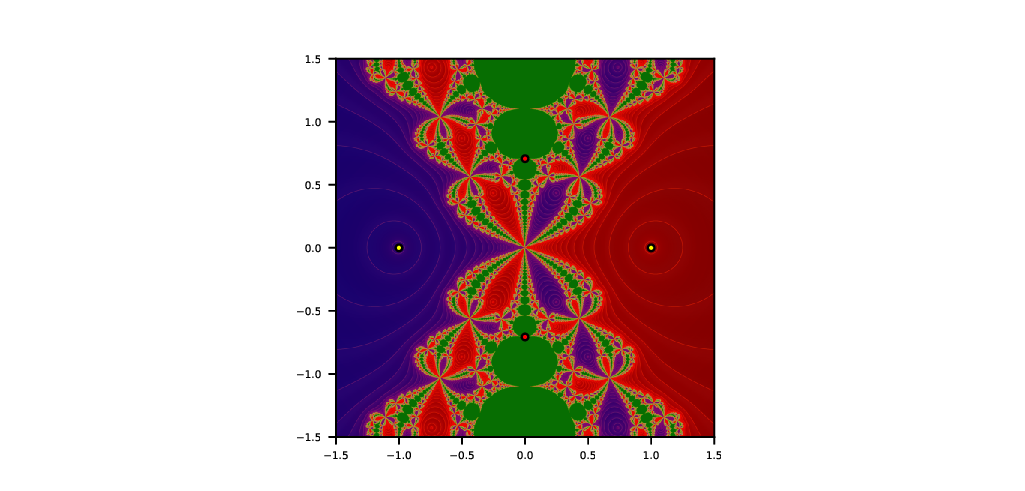}
    \caption{$T_{p_{\vphantom{1}_{2}},-\frac{4}{3}}$}
    \label{fig:T_p_2,-4/3}
  \end{subfigure}
  \caption{Dynamical plane of $T_{p_{\vphantom{1}_{2}},\delta}$ when $\delta\in(-4,-\frac{4}{3}]$. The blue region is $\mathcal A(-1)$ and the red region is $\mathcal A(1)$. The red dots are the two extraneous fixed points $e_1$ and $e_2$ and the green regions indicate their basins.}
  \label{fig:extra}
\end{figure} 
\begin{proof}[Proof of Theorem \ref{T:(-4/3,0)}]
  %  As $\delta\in(-\frac{4}{3},0)$, for all $x>1$, $T_{p_{\vphantom{1}_{2}},\delta}(x)<x$ and $T_{p_{\vphantom{1}_{2}},\delta}^\prime(x)>0$. Also, $T_{p_{\vphantom{1}_{2}},\delta}(1)=1$. So, $T_{p_{\vphantom{1}_{2}},\delta}^n(x)\to 1$ for all $x\geq 1$. Therefore, $[1,\infty)\subset \mathcal A_1$. Hence, $(-\infty,-1]\subset \mathcal A_{-1}$. So, $\mathcal A_{-1}$ and $\mathcal A_{1}$ are invariant unbounded immediate attracting basins of $T_{p_{\vphantom{1}_{2}},\delta}$. Since $\delta\in(-\frac{4}{3},0)$, the point at infinity is a repelling fixed point. So, $0\in J(T_{p_{\vphantom{1}_{2}},\delta})$. For $\delta\in(-\frac{4}{3},0)$, both free critical points are purely imaginary. Therefore by Lemma \ref{L:Tp2invline}, they do not lie in either $\mathcal A_{-1}$ or $\mathcal A_1$. So, both $\mathcal A_{-1}$ and $\mathcal A_1$ are simply connected by \cite[Proposition 3.1]{campos2020connectivity}.
   % Therefore, by \cite[Lemma 3.4]{nayak2025symmetry}, $0\in \partial \mathcal A_1$ and $J(T_{p_{\vphantom{1}_{2}},\delta})$ is connected.

For $\delta\in\left(-\frac{4}{3},0\right)$, we have, for every $x>1$,
\[
T_{p_2,\delta}(x)<x
\quad\text{and}\quad
T_{p_2,\delta}'(x)>0.
\]
 It follows from Eqn. \eqref{E:Tp2inverse-1} that $T_{p_2,\delta}(1)=1$. Since $T_{p_2,\delta}$ is strictly increasing on $(1,\infty)$, we get
\[
T_{p_2,\delta}(x)>T_{p_2,\delta}(1)=1
\]
for every $x>1$ concluding 
\[
1<T_{p_2,\delta}(x)<x,\qquad x>1.
\]
Consequently, for every $x>1$, the sequence of iterates $\{T_{p_2,\delta}^{n}(x)\}_{n\geq 0}$ is decreasing and bounded below by $1$. Hence, it converges to some $\ell\geq 1$. By continuity of $T_{p_2,\delta}$, we have $T_{p_2,\delta}(\ell)=\ell$. Since $T_{p_2,\delta}(x)<x$ for all $x>1$, the only possible limit is $\ell=1$. Therefore,
\[
T_{p_2,\delta}^{n}(x)\longrightarrow 1
\quad\text{as }n\to\infty
\]
for every $x\geq 1$. It follows that $[1,\infty)\subset\mathcal A_1$. By the symmetry of $T_{p_2,\delta}$, we similarly obtain $(-\infty,-1]\subset\mathcal A_{-1}$. In particular, both attracting basins are unbounded, and their immediate components containing $1$ and $-1$, respectively, are unbounded.

Moreover, for $\delta\in\left(-\frac{4}{3},0\right)$, the point at infinity is a repelling fixed point of $T_{p_2,\delta}$, and hence, $\infty\in J(T_{p_2,\delta})$. By the symmetry of the Julia set, it follows that $0\in J(T_{p_2,\delta})$. For $\delta\in\left(-\frac{4}{3},0\right)$, both free critical points of $T_{p_2,\delta}$ are purely imaginary. Hence, by Lemma \ref{L:Tp2invline}, neither of them belongs to $\mathcal A_{-1}$ or $\mathcal A_1$. Consequently, by \cite[Proposition 3.1]{campos2020connectivity}, both immediate attracting basins $\mathcal A_{-1}$ and $\mathcal A_1$ are simply connected. Finally, by \cite[Lemma 3.4]{nayak2025symmetry}, we conclude that $0\in\partial\mathcal A_1$ and that the Julia set $J(T_{p_2,\delta})$ is connected.

\end{proof}
\begin{proof}[Proof of Theorem \ref{T:0,4}]
For $\delta\in(0,4)$, the point at infinity lies in the Julia set and hence, $0\in J(T_{p_{\vphantom{1}_{2}},\delta})$.
\begin{enumerate}[label=\alph*)]

\item 

Let $I=\{iy:y\in\mathbb R\}$ denote the imaginary axis. We claim that
\[
I\subset J(T_{p_2,\delta}).
\]
Suppose, to the contrary, that $I\not\subset J(T_{p_2,\delta})$. Then there exists a Fatou component $U$ such that $U\cap I\neq\emptyset$.

Since $T_{p_2,\delta}$ is a rational map and its Fatou set contains neither Siegel disks nor Herman rings, every Fatou component eventually maps onto a periodic immediate attracting basin or a periodic immediate parabolic basin. Thus, there exists an integer $k\geq 0$ such that $T_{p_2,\delta}^{k}(U)=V$, where $V$ is a $p$-periodic immediate attracting basin or an immediate parabolic basin. Moreover, at least one of the free critical points $c_1$ or $c_2$ belongs to the cycle of $V$. Without loss of generality, assume it belongs to $V$. By Lemma \ref{L:Tp2invline}, we have
\[
V\cap I\neq\emptyset.
\]

 Since $V\cap I\neq\emptyset$, there exists $z_0=iy\in V$. Observe that $-\overline{z_0}=iy=z_0$.
Define $V'=\{-\overline z:z\in V\}$. Hence, $z_0\in V\cap V'$, which implies that $V\cap V'\neq\emptyset$ concluding that $V=V'$. In particular, $V$ is invariant under the reflection $z\mapsto-\overline z$. For $\delta\in(0,4)$, both free critical points $c_1$ and $c_2$ are real, and $c_2=-c_1$. As at least one of the two critical points belongs to $V$, both $c_1$ and $c_2$ belong to $V$. Since $V$ is a $p$-periodic immediate attracting basin or an immediate parabolic basin, there exists a point $x_0\in\mathbb R$ such that
\[
T_{p_2,\delta}^{np}(c_1)\longrightarrow x_0
\quad\text{as }n\to\infty.
\]
By the symmetry of $T_{p_2,\delta}$, we also have
\[
T_{p_2,\delta}^{np}(c_2)
=
-T_{p_2,\delta}^{np}(c_1).
\]
Therefore,
\[
T_{p_2,\delta}^{np}(c_2)\longrightarrow -x_0.
\]
On the other hand, since $c_2\in V$, its iterates under $T_{p_2,\delta}^{p}$ converge to the same periodic limit point $x_0$. Thus, $x_0=-x_0$ giving that $x_0=0$. This contradicts the fact that $0\in J(T_{p_2,\delta})$ is not a periodic point. Hence, our assumption was false, and we conclude that $I\subset J(T_{p_2,\delta})$.

    \item For $\delta\in(0,1)$, we have $0<\frac{3\delta}{4-\delta}<1$. From Eqn. \eqref{E:Tp2inverse1}, $T_{p_{\vphantom{1}_{2}},\delta}(c_1)>1$. For all $x>1$, $T_{p_{\vphantom{1}_{2}},\delta}(x)-x<0$ and $T_{p_{\vphantom{1}_{2}},\delta}^\prime(x)>0$. As $T_{p_{\vphantom{1}_{2}},\delta}(1)=1$, $T_{p_{\vphantom{1}_{2}},\delta}^n(x)\to 1$ for all $x\geq 1$. Hence, $T_{p_{\vphantom{1}_{2}},\delta}^n(c_1)\to 1$ as $n\to \infty$. Using Lemma \ref{L:Tp2symmetry}, we have $T_{p_{\vphantom{1}_{2}},\delta}^n(c_2)\to -1$ as $n\to \infty$. Therefore, $c_1\in \mathcal A(1)$ and $c_2\in \mathcal A(-1)$. As $\delta\in (0,1)$, $\infty\in J(T_{p_{\vphantom{1}_{2}},\delta})$ and consequently,  $0\in J(T_{p_{\vphantom{1}_{2}},\delta})$. Hence $$F(T_{p_{\vphantom{1}_{2}},\delta})=\mathcal A(-1)\cup \mathcal A(1).$$ As $T_{p_{\vphantom{1}_{2}},\delta}^n(x)\to 1$ for all $x\geq 1$, $[1,\infty)\subset \mathcal A(1)$. Since, $[1,\infty)]$ is connected, $[1,\infty)\subset \mathcal A_1$. Hence, $\mathcal A_1$ is unbounded. Consequently, $\mathcal A_{-1}$ is also unbounded.    
    For $\delta\in(0,1)$, $T_{p_{\vphantom{1}_{2}},\delta}^\prime(x)>0$ for all $x\in(c_2,0)$ and $T_{p_{\vphantom{1}_{2}},\delta}^\prime(x)<0$ for all $x\in(-1,c_2)$. The map $T_{p_{\vphantom{1}_{2}},\delta}$ has a local minimum at $c_2$. By Lemma \ref{T:Tp2civ}, we have $T_{p_{\vphantom{1}_{2}},\delta}(r_{-1})=-1$ and $r_{-1}>c_2$. Thus for all $x\in (-1, r_{-1})$, $T_{p_{\vphantom{1}_{2}},\delta}(x)<-1$. This gives that $(-1, r_{-1})\subset \mathcal A(-1)$. The point $r_{-1}\in \mathcal A(-1)$ implies $[-1,r_{-1}]\subset \mathcal A(-1)$. As $[-1,r_{-1}]$ is connected, we have $[-1,r_{-1}]\subset \mathcal A_{-1}$. The Fatou component $\mathcal A_{-1}$ contains a critical point $c_2(\neq -1)$ and a preimage of $-1$, $r_{-1}(\neq -1)$. Hence, by \cite[Proposition 3.1]{campos2020connectivity}, the Fatou component $\mathcal A_{-1}$ is simply connected. If $U\neq \mathcal A_{-1}$ is a component of $T_{p_{\vphantom{1}_{2}},\delta}^{-1}(\mathcal A_{-1})$ then $U$ contains no critical points concluding that   $U$ is simply connected (\cite[Corollary 2.2]{canela2015family}). Any preperiodic Fatou component landing on $\mathcal A_{-1}$ has no critical points, all components of $\mathcal A(-1)$ are simply connected. Similarly, all components of $\mathcal A(1)$ are simply connected. As the Fatou set of $T_{p_{\vphantom{1}_{2}},\delta}$ is the union of $\mathcal A(-1)$ and $\mathcal A(1)$, all Fatou components are simply connected. Hence, $J(T_{p_{\vphantom{1}_{2}},\delta})$ is connected for all $\delta\in (0,1)$. 
    
    For $\delta=1$, the only critical points of $T_{p_{\vphantom{1}_{2}},1}$ are $-1,1$ and $0$. As $0\in J(T_{p_{\vphantom{1}_{2}},1})$, $$F(T_{p_{\vphantom{1}_{2}},\delta})=\mathcal A(-1)\cup \mathcal A(1).$$ 
    In this case, $\mathcal A_{-1}$ contains no critical point other than $-1$ and $\mathcal A_1$ contains no critical point other than $1$. So, by \cite[Proposition 3.1]{campos2020connectivity}, they are simply connected. The rest of the proof follows as in case of $\delta\in(0,1)$.
    \item For $\delta\in(1,\frac{3+\sqrt{13}}{2})$, we have $T_{p_{\vphantom{1}_{2}},\delta}(c_1)=\frac{(2\delta+1)\sqrt{4-\delta}}{3\sqrt{3\delta}}$ and $e_1<T_{p_{\vphantom{1}_{2}},\delta}(c_1)<1$.  Note that $T_{p_{\vphantom{1}_{2}},\delta}(1)=1$. Moreover, $T_{p_{\vphantom{1}_{2}},\delta}(x)>x$ and $T_{p_{\vphantom{1}_{2}},\delta}^\prime(x)>0$  for all $x\in(e_1,1)$ concluding that $T_{p_{\vphantom{1}_{2}},\delta}^n(x)\to 1$ for all $x\in (e_1,1]$. The critical value $T_{p_{\vphantom{1}_{2}},\delta}(c_1)\in \mathcal A_1$ and $c_1\in \mathcal A(1)$,  and hence,  $c_2\in \mathcal A(-1)$. As all the critical points are in  $\mathcal A(1)\cup \mathcal A(-1)$, 
    $$F(T_{p_{\vphantom{1}_{2}},\delta})=\mathcal A(1)\cup \mathcal A(-1).$$ 
  From Remark \ref{R:1pre}, we have $r_1>c_1>1$. The map $T_{p_{\vphantom{1}_{2}},\delta}$ attains a local minimum at $c_1$ in $(1,r_1)$. Further, $T_{p_{\vphantom{1}_{2}},\delta}$ is decreasing in $(1,c_1)$ and increasing in $(c_1,r_1)$. Since $T_{p_{\vphantom{1}_{2}},\delta}(c_1)\in (e_1,1)$, we have $T_{p_{\vphantom{1}_{2}},\delta}(1,r_1)\subset (e_1,1)$ and hence, $(1,r_1]\subset \mathcal A(1)$. Note that  $(1,r_1]\subset \mathcal A_1$ as $(1,r_1]$ is connected. The Fatou component $\mathcal A_1$ contains the critical point $c_1\neq 1$ and a preimage $r_1$ of $1$ giving that $\mathcal A_1$ is simply connected. Similarly, $\mathcal A_{-1}$ is also simply connected. The rest of the proof follows as in the last case.\par
    For $\delta=\frac{3+\sqrt{13}}{2}$, we have $T_{p_{\vphantom{1}_{2}},\delta}(c_1)=e_1$. Therefore $\mathcal A_1$ and $\mathcal A_{-1}$ contains exactly one critical point each concluding that both $\mathcal A_1$ and $\mathcal A_{-1}$ are simply connected. Now, the proof follows as before.
\end{enumerate}     
\end{proof}
\begin{figure}[H]
  \centering  
  \begin{subfigure}[b]{0.48\textwidth}
    \centering
    \includegraphics[trim = 5cm 0cm 5cm 0cm, clip, scale=1]{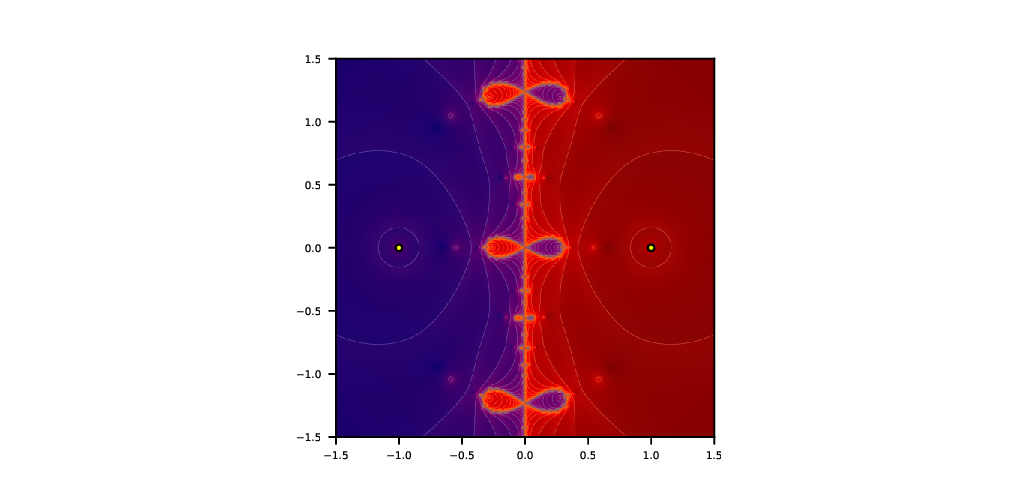}
    \caption{$T_{p_{\vphantom{1}_{2}},0.5}$}
    \label{fig:T_p_2,0.5}
  \end{subfigure}
  \hfill
  \begin{subfigure}[b]{0.48\textwidth}
    \centering
    \includegraphics[trim = 5cm 0cm 5cm 0cm, clip, scale=1]{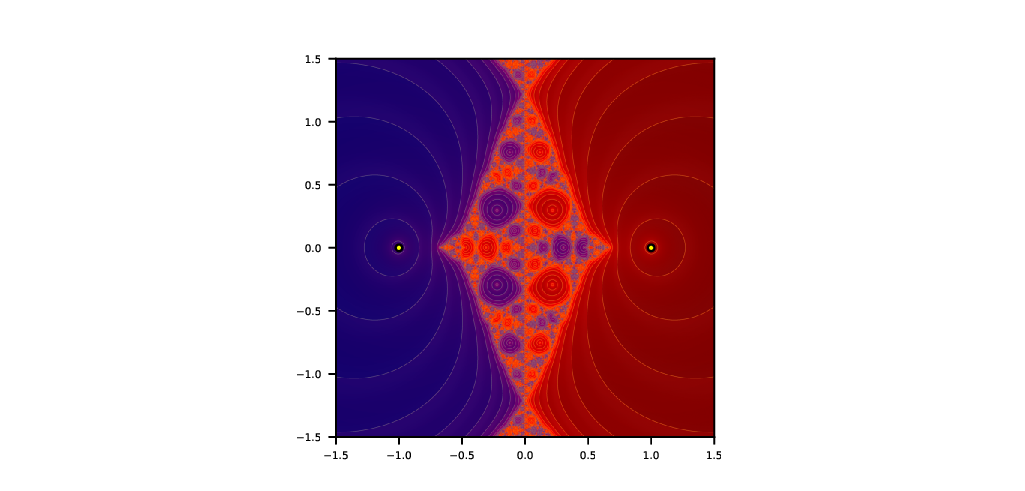}
    \caption{$T_{p_{\vphantom{1}_{2}},3.302775}$}
    \label{fig:T_p_2,3+root13/2}
  \end{subfigure}
  \caption{Dynamical plane of $T_{p_{\vphantom{1}_{2}},\delta}$ when $\delta\in(0,4)$. The blue region is $\mathcal A(-1)$ and the red region is $\mathcal A(1)$.}
  \label{fig:0.5}
\end{figure}
 We have observed that for $\delta\in (\frac{3+\sqrt{13}}{2},4)$, $c_1$ can converge to $-1$ while $c_2$ can converge to $1$. There are also some values of $\delta$ for which the free critical orbits enter higher periodic cycles. However, our method is unable to show this.

Numerical observations show that the free critical orbits enter higher periodic cycles for $\delta\in(4,\frac{19+5\sqrt{17}}{8})$. Due to the complicated nature of the higher iterates of $T_{p_{\vphantom{1}_{2}},\delta}$, analytically showing this might not be possible. 

\begin{proof}[Proof of Theorem \ref{T:4.8,12}]
For $\delta\in[\frac{19+5\sqrt{17}}{8},12)$, both free critical points are purely imaginary.  We need to analyze the map $g(y)$ described in Lemma \ref{L:Tp2invline}. It is obvious from the definition that $g$ has no real fixed points. Assuming $g(y)=-y$, we have $g^2(y)=g(-y)=-g(y)=y$ giving that $g$ has a two-periodic point at $y$. Solving the equation $g(y)=-y$, we get  $y=\pm\sqrt{\frac{2+\delta+2\sqrt{4\delta+1}}{12-\delta}}\in\mathbb{R}$. The $2$-periodic points of $g$ are  $$y_1=\sqrt{\frac{2+\delta+2\sqrt{4\delta+1}}{12-\delta}} \text{ and } -y_1=-\sqrt{\frac{2+\delta+2\sqrt{4\delta+1}}{12-\delta}}.$$  Upon calculation, we have $g^\prime(y_1)=\frac{7\delta+2-(2+\delta)\sqrt{4\delta+1}}{\delta}$ and $g^{2^\prime}(y_1)=g^\prime(g(y_1))g^\prime(y_1)=(g^\prime(y_1))^2$. Considering $$M(\delta)=\frac{7\delta+2-(2+\delta)\sqrt{4\delta+1}}{\delta},$$ the derivative is $$M'(\delta)=-2\frac{(\delta-1)^2+\sqrt{4\delta+1}-2}{\delta^2\sqrt{4\delta+1}}<0$$ for $4<\delta<12$. Note that  $M(\frac{19+5\sqrt{17}}{8})=1$ and $M(12)=-1$.
\begin{enumerate}[label=\alph*)]
    \item For $\delta\in(\frac{19+5\sqrt{17}}{8},12)$, we get $1>M(\delta)>-1$ as $M$ is a decreasing giving that $\lvert g^{2^\prime}(y_1)\rvert<1$. Clearly, $T_{p_{\vphantom{1}_{2}},\delta}^2(iy_1)=iy_1$ and $T_{p_{\vphantom{1}_{2}},\delta}^{2^\prime}(iy_1)=g^{2^\prime}(y_1)$ concluding that $iy_1$ is an attracting $2$-periodic point of $T_{p_{\vphantom{1}_{2}},\delta}$. Similarly, $-iy_1$ is also an attracting $2$-periodic point of $T_{p_{\vphantom{1}_{2}},\delta}$. Thus, there exists $2$-cycle attracting basins $\{U,-U\}$ such that $iy_1\in U$. Since $U$ is an attracting basin, $U$ must contain one free critical point and hence, $-U$ contains the other free critical point. Hence,     $$F(T_{p_{\vphantom{1}_{2}},\delta})=\mathcal A(-1)\cup \mathcal A(1)\cup U\cup (-U).$$
    
    \item The next part follows similarly. In this case, $M(\frac{19+5\sqrt{17}}{8})=1$ and hence, $iy_1$ and $-iy_1$ are rationally indifferent periodic points concluding that  there exists $2$-cycle parabolic basins $\{V,-V\}$. Similar to the above, we have     $$F(T_{p_{\vphantom{1}_{2}},\delta})=\mathcal A(-1)\cup \mathcal A(1)\cup V\cup (-V).$$
    
    \item For $\delta\in(\frac{19+5\sqrt{17}}{8},12)$, let $iy_1\in U^1$ where $U^1$ is Fatou component. So, $U^1\subset U$ and $-U^1$ is the Fatou component containing $-iy_1$. Clearly, $U^1\neq -U^1$. Otherwise $U^1$ contains two distinct periodic points $iy_1$ and $-iy_1$ which is not possible. Since, $\{U^1,-U^1\}$ is a $2$-cycle of immediate attracting basins, at least one free critical point lie in $U^1$ or $-U^1$. Since $U^1\neq -U^1$, both $U^1$ and $-U^1$ contains exactly one simple critical point due to Lemma \ref{L:Tp2symmetry}. Let $D\subset U^1$ be a sufficiently small disk containing $iy_1$ such that $T_{p_{\vphantom{1}_{2}},\delta}^2(D)\subset D$ and $T_{p_{\vphantom{1}_{2}},\delta}(D)$ is a disk. This is possible as $iy_1$ is attracting for all $\delta\in(\frac{19+5\sqrt{17}}{8},12)$. Let $A^0=D$ and $B^0=T_{p_{\vphantom{1}_{2}},\delta}(D)$. Also, let $A^n$ be the component of $T_{p_{\vphantom{1}_{2}},\delta}^{-1}(B^{n-1})$ containing $A^{n-1}$ and $B^n$ be the component of $T_{p_{\vphantom{1}_{2}},\delta}^{-1}(A^{n-1})$ containing $B^{n-1}$. So, $A^n\subset U^1$ and $B^n\subset -U^1$. So, each $A^n$ or $B^n$ has at most one critical point. Therefore by induction and \cite[Corollary 2.2]{canela2015family}, each $A^n$ and $B^n$ is simply connected. Let $z\in U^1$, $w\in D$ and $\gamma$ be a path joining $z$ and $w$. Then there exists $k\in\mathbb{N}$ such that $T_{p_{\vphantom{1}_{2}},\delta}^{2k}(\gamma)\subset D$. This implies that for some $n\in\mathbb{N}$, $z\in A^n$. Therefore, $\displaystyle U^1=\bigcup_{n\in\mathbb{N}}A^n$ and $\displaystyle -U^1=\bigcup_{n\in\mathbb{N}}B^n$. Hence, $U^1$ and $-U^1$ are both simply connected and $\mathcal A_1,\mathcal A_{-1}$ are also simply connected by \cite{campos2020connectivity}. This gives that all the Fatou components are simply connected using Riemann-Hurwitz formula. Hence, the Julia set is connected. The proof follows similarly for $\delta=\frac{19+5\sqrt{17}}{8}$.
\end{enumerate}  
\end{proof}
\begin{figure}[H]
    \centering
    \includegraphics[trim = 5cm 0cm 5cm 0cm, clip, scale=1]{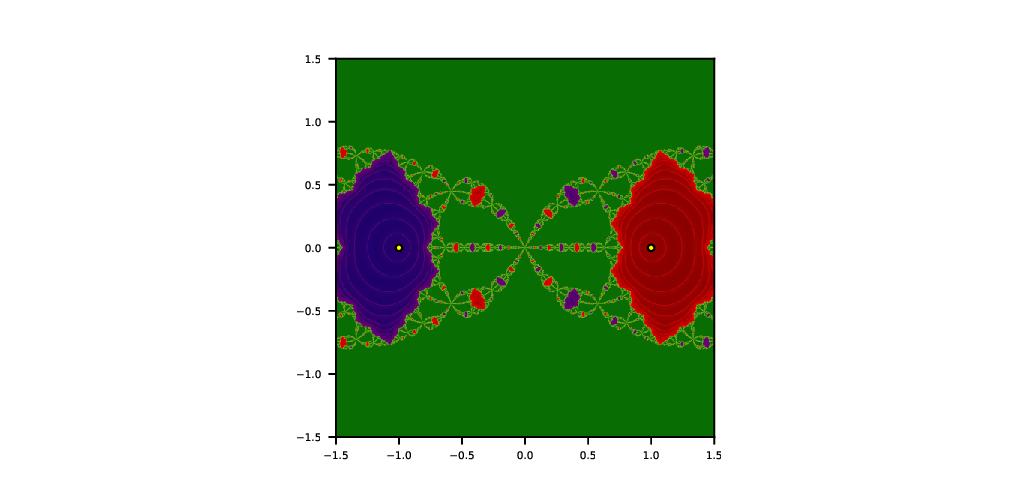}
    \caption{Dynamical plane of $T_{p_{\vphantom{1}_{2}},6}$. The blue region is $\mathcal A(-1)$ and the red region is $\mathcal A(1)$. The green regions denote the basins of the attracting $2$-periodic points.}
    \label{fig:T_p_2,6}
\end{figure}
As mentioned before, for $\delta\in(4,\frac{19+5\sqrt{17}}{8})$ the free critical points converge to higher periodic cycles. So, we use Python programming to give an estimate of these cycles for different values of $\delta\in(4,\frac{19+5\sqrt{17}}{8})$ in Table \ref{tab:1}.
\newpage
    \begin{table}[h!]
    \centering
    \captionsetup{width=\textwidth}
    \renewcommand{\arraystretch}{1.8}
    \begin{tabular}{|c|c|m{6cm}|}
    \hline
        $\delta$ & Period & \begin{center}Cycles of periodic points of $T_{p_{\vphantom{1}_{2}},\delta}$\end{center}\\
        \hline
         $4.8$ & $2$ & \makecell{$- 0.9546i$, $2.5658i$} \\
         \hline
         $4.1$ & $4$ & \makecell{$0.2770i$,  $-29.6209i$, $0.4218i$,\\ $-10.4519i$}\\
         \hline
         $4.07$ & 8 & \makecell{$-0.2493i$, $38.9179i$, $- 0.3795i$,\\ $13.3074i$, $- 0.2306i$, $48.0199i$\\
$- 0.4517i$, $8.8801i$}\\
         \hline
         $4.06$ & 16 & \makecell{$- 0.2465i$, $40.0016i$,
$- 0.3378i$,\\ $17.6413i$,
$- 0.2182i$, $55.7392i$,\\
$- 0.4452i$, $9.1564i$,
$- 0.2347i$,\\ $45.6632i$,
$- 0.3756i$, $13.6090i$,\\
$- 0.2135i$, $59.1755i$,
$- 0.4694i$,\\ $8.1371i$}\\
\hline
    \end{tabular}
    \caption{Periodic point cycles of $T_{p_{\vphantom{1}_{2}},\delta}$ for $\delta\in(4,\frac{19+5\sqrt{17}}{8})$.}
    \label{tab:1}
\end{table}
\begin{proof}[Proof of Theorem \ref{T:12,inf}]
If $\delta\in(12,\infty)$, the point at infinity is an attracting fixed point of $T_{p_{\vphantom{1}_{2}},\delta}$ and for $\delta=12$, it is a rationally indifferent fixed point by Remark \ref{R:rangeinfty}.
\begin{enumerate}[label=\alph*)]
    \item As $\delta\in(12,\infty)$, we have $0\in \mathcal A(\infty)$. In this case, both free critical points are purely imaginary. From Lemma \ref{L:Tp2invline} we have, 
    $g(y)=y\left(1-\frac{(y^2+1)\big((4+\delta)y^2+\delta\big)}{8y^4}\right).$
    Hence, 
    $$\frac{g(y)}{y}=-\Big( \big(\frac{\delta}{8}-\frac{1}{2}\big)+\frac{2+\delta}{4y^2}+\frac{\delta}{8y^4}\Big).$$ 
     For $y>0$, we have $g(y)<0$ and for $y<0$, the function $g(y)>0$. Clearly,  $\frac{\delta}{8}-\frac{1}{2}>1$ as $\delta>12$. This gives that  $\Big\lvert\frac{g(y)}{y}\Big\rvert>1$ for all $y\neq 0$. In other words, $\lvert g(y)\rvert>\lvert y\rvert$ for all $y\neq 0$. So, $\{\lvert g^n(y)\rvert\}$ is an unbounded strictly increasing sequence for all $y$ such that $g^n(y)\neq 0$ for all $n\in\mathbb{N}$. Thus, $g^n(y)$ keeps oscillating between the negative real numbers and positive real numbers while steadily increasing in magnitude for all such $y$.  Therefore, $T_{p_{\vphantom{1}_{2}},\delta}^n(iy)\to \infty$ as $n\to\infty$ for all $y\neq 0$ concluding that both free critical points lie in $\mathcal A_\infty$. Thus, all the critical points lie in one of $\mathcal A(-1), \mathcal A(1)$ or $\mathcal A_\infty$ and hence,  $$F(T_{p_{\vphantom{1}_{2}},\delta})=\mathcal A(-1)\cup \mathcal A(1)\cup \mathcal A(\infty).$$ 

     It follows from above discussion that the set $\{iy:y\neq 0\}$ is contained in  $\mathcal A_\infty$ for $\delta\in (12,\infty)$. Also, $T_{p_{\vphantom{1}_{2}},\delta}(0)=\infty$ giving that $0\in F(T_{p_{\vphantom{1}_{2}},\delta})$ lies in $\mathcal A(\infty)$. This implies $\{0\}\cup \{iy:y\neq 0\}\subset \mathcal A_\infty$ i.e., the imaginary axis is contained in $\mathcal A_\infty$. Hence, the Julia set of  $T_{p_{\vphantom{1}_{2}},\delta}$ is disconnected. 
    
    \item For $\delta=12$,  the proof of the first part is similar to that of Theorem \ref{T:12,inf} (a). For the second part, we use the semiconjugacy relation defined in Proposition \ref{P:semiconj}. For $\delta=12$, the map $R(z)=\frac{(2z^2-7z+3)^2}{4z^3}$. The map $R$ has critical points at $0,1,\frac{1}{2},3$, and $-\frac{9}{2}$. Now, $0$ is a pole of $R$ and $R(\frac{1}{2})=R(3)=0$. So, $0,\frac{1}{2},3$ all lie in the Julia set. Also, $R(1)=1$ and $R^\prime(1)=0$. So, $1$ is a superattracting fixed point of $R$. For the map $R$, $\infty$ is a rationally indifferent fixed point with multiplier $1$. So, the immediate parabolic basin of $\infty$, $V$ must contain $-\frac{9}{2}$. Let $P$ be a simply connected invariant petal in $V$ such that $\infty\in\partial P$, and all points in $V$ eventually enters $P$ under iteration. Let $P^0=P$ and $P^{n}\subset V$ be the component of $R^{-1}(P^{n-1})$ that contains $P$. Since, $P^n$ contains at most one critical point $-\frac{9}{2}$ for any $n\in\mathbb{N}\cup\{0\}$, by \cite[Corollary 2.2]{canela2015family} and mathematical induction each $P^n$ is simply connected. Let $z\in V$ and $\gamma\subset V$ joining $z$ to a point of $P$. Then there exists some $n$ such that $R^n(\gamma)\subset P$. So, $\gamma\cup P$ is a connected subset of $R^{-n}(P)$ containing $P$. So, $z\in P^n$. This implies     
    $\displaystyle V=\bigcup_{n\in\mathbb{N}\cup\{0\}} P^n$. So, $V$ is union of nested simply connected domains. Thus, $V$ is simply connected and $\partial V$ is connected. Now, $\partial V$ contains $\infty$. So, by Lemma \ref{L:general1}, $0\in\partial V$. Consequently by Lemma \ref{L:general2}, $J(R)$ is connected. Now, $J(T_{p_{\vphantom{1}_{2}},12})=\phi^{-1}(J(R))$. 
    $J(T_{p_{\vphantom{1}_{2}},12})$ is compact, contains $0$ and is symmetric by Lemma \ref{L:Tp2symmetry}. Let $J(T_{p_{\vphantom{1}_{2}},12})$ be disconnected. This implies that $J(T_{p_{\vphantom{1}_{2}},12})=A\cup B$ where $A$ and $B$ are nonempty, disjoint and closed sets. Suppose $0\in A$. Then $A\cap (-A)$ is symmetric closed set and contains $0$. Also, $J(T_{p_{\vphantom{1}_{2}},12})\setminus (A\cap(-A))$ is symmetric and closed. Now, if $\phi(z)=\phi(w)$ such that $z\in A\cap(-A)$ and $w\in J(T_{p_{\vphantom{1}_{2}},12})\setminus (A\cap(-A))$, then we have $z=\pm w$ which implies $w\in A\cap(-A)$. This is not possible. So, $\phi(A\cap(-A))$ and $\phi(J(T_{p_{\vphantom{1}_{2}},12})\setminus (A\cap(-A)))$ are disjoint and $J(R)= \phi(A\cap(-A))\cup \phi(J(T_{p_{\vphantom{1}_{2}},12})\setminus (A\cap(-A)))$. Since, $J(R)$ is connected. This implies one of $\phi(A\cap(-A))$ or $\phi(J(T_{p_{\vphantom{1}_{2}},12})\setminus (A\cap(-A)))$ is empty. So, $A\cap (-A)=J(T_{p_{\vphantom{1}_{2}},12})$. Then $A=J(T_{p_{\vphantom{1}_{2}},12})$ and $B=\emptyset$. This contradicts $B$ is non-empty. This implies $J(T_{p_{\vphantom{1}_{2}},12})$ is connected. 
    \end{enumerate}    
\end{proof}
\section*{Statements and Declarations}

\subsection*{Competing Interests} The authors have no competing interests to declare that are relevant to the content of this article.

\subsection*{Data availability statement} Not applicable to this article as no datasets were generated or analyzed during the work.

\subsection*{Funding}
S.C. and G.C. would like to acknowledge the financial support rendered by DST-FIST, Govt. of India (File no: SR/FST/MSI/2022/123) running at the Department of Mathematics, Sidho-Kanho-Birsha University, Purulia, 723104, India. 
\bibliographystyle{plain}
\bibliography{sn-bibliography}
\end{document}